\documentclass[12pt]{article}
\usepackage{diagbox}
\usepackage{mathtools}
\usepackage{bbm}
\usepackage{amsmath}
\usepackage{hyperref}
\allowdisplaybreaks

\usepackage[utf8x]{inputenc}

\renewcommand{\emph}[1]{\textit{#1}}
 \usepackage{enumerate,amsmath,amsthm,latexsym,amssymb}
\usepackage{color}\usepackage{graphicx}
\newcommand{\bi}{\boldsymbol\sigma}
\def\bM{{\bf M}}
\def\bm{\bf m}

\definecolor{brown}{cmyk}{0, 0.72, 1, 0.45}
\definecolor{grey}{gray}{0.5}

\newcommand{\old}[1]{}

\newcounter{rot}

\newcommand{\card}[1]{\left|#1\right|}

\newcommand{\ignore}[1]{}

\def\bM{{\bf M}}

\def\cA{{\mathcal A}}

\newcommand{\set}[1]{\left\{#1\right\}}

\def\cP{\mathcal{P}}

\def\ii_(#1,#2){i_{#1}^{#2}}
\newcommand{\cdnm}[1]{{\mathcal D}_{n,m}^{(\delta\geq#1)}}
\newcommand{\dnm}[1]{{D}_{n,m}^{(\delta\geq#1)}}

\def\a{\alpha}
\def\b{\beta}
\newcommand{\de}{\delta}
\newcommand{\D}{\Delta}
\def\e{\varepsilon}
\def\f{\phi}

\def\g{\gamma}
\def\G{\Gamma}
\def\k{\kappa}

\def\th{\theta}

\def\l{\lambda}
\def\m{\mu}
\def\n{\nu}
\def\p{\pi}

\def\r{\rho}

\def\s{\sigma}
\def\t{\tau}
\def\om{\omega}

\def\w{\omega}

\def\1{{\bf 1}}
\def\0{{\bf 0}}

\newcommand{\rdup}[1]{\left\lceil #1 \right\rceil}
\newcommand{\rdown}[1]{\mbox{$\left\lfloor #1 \right\rfloor$}}

\newcommand{\wh}[1]{\widehat #1}

\def\cE{\mathcal{E}}
\def\cW{\mathcal{W}}
\def\cF{\mathcal{F}}
\def\cL{\mathcal{L}}

\newcommand{\brac}[1]{\left( #1 \right)}

\def\E{{\bf E}}
\def\Var{{\bf Var}}
\renewcommand{\Pr}{\operatorname{\bf Pr}}
\newcommand\bfrac[2]{\left(\frac{#1}{#2}\right)}

\newcounter{Prob}
\newtheorem{theorem}{Theorem}[section]

\newtheorem{remthm}[theorem]{Remark}

\newcounter{thmtemp}
\newtheorem{Problem}[Prob]{Problem}

\newcommand{\nospace}[1]{}

\def\path{\operatorname{PATH}}

\newcommand{\beq}[2]{\begin{equation}\label{#1}#2\end{equation}}

\def\va{\boldsymbol \a}

\def\bc{{\bf c}}
\def\cG{{\cal G}}

\def\bd{{\bf d}}

\def\cX{{\cal X}}

\def\cH{{\cal H}}
\usepackage{tikz}
\usetikzlibrary{arrows}
\usetikzlibrary{decorations}
\usetikzlibrary{shapes.misc}
\newcommand{\pic}[1]{
\begin{tikzpicture}
#1
\end{tikzpicture}}

\author{Alan Frieze\thanks{Research supported in part by NSF grant DMS1661063
} \\Department of Mathematical Sciences\\Carnegie Mellon University\\Pittsburgh PA 15213} 
\begin{document}
\title{Hamilton Cycles in Random Graphs: a bibliography}
\maketitle
\begin{abstract}
We provide an annotated bibliography for the study of Hamilton cycles in random graphs and hypergraphs.
\end{abstract}
\section{Introduction}
As is well-known, the study of the structure of random graphs began in earnest with two seminal papers by Erd\H{o}s and R\'enyi \cite{ER1}, \cite{ER2}. At the end of the \cite{ER2} the authors pose the question: ``for what order of magnitude of $N(n)$ has $\G_{n,N(n)}$ with probability tending to 1 a {\em Hamilton-line} (i.e. a path which passes through all vertices)''. Thus began the study of Hamilton cycles in random graphs. By now there is an extensive literature on this and related problems and the aim of this paper to summarise what we know and what we would like to know about these questions.

{\bf Notation:} Our notation for random graphs is standard and can be found in any of Bollob\'as \cite{Boll1}, Frieze and Karo\'nski \cite{FK} or Janson, {\L}uczak and Rucinski \cite{JLR}.
\section{The random graphs $G_{n,m}$ and $G_{n,p}$}
\subsection{Existence}
In this section we consider the random graphs $G_{n,m}$, $G_{n,p}$ and the random process $G_m,m=0,1,\ldots,N=\binom{n}{2}$. The first paper to make significant progress on the threshold for Hamilton cycles was by Koml\'os and Szemer\'edi \cite{KS1} who proved that $m=n^{1+\e}$ is sufficent for any positive constant $\e>0$. A breakthrough came when Pos\'a \cite{Posa} showed that $m=O(n\log n)$ is sufficient and introduced the idea of using rotations. Given a longest path $P=(x_1,x_2,\ldots,x_s)$ in a graph $G$ and an edge $\set{x_s,x_j},1<j<s-1$ we can create another longest path $P'=(x_1,x_2,\ldots,x_j,x_s,x_{s-1},\ldots,x_{j+1})$ with a new endpoint $x_{j+1}$. We call this a rotation.
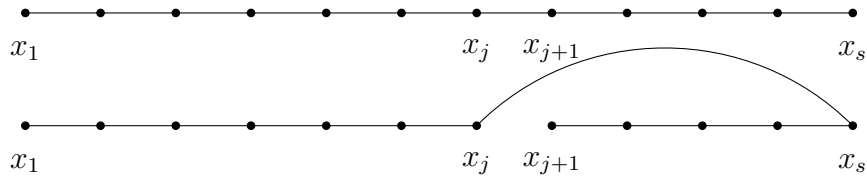
\begin{center}
\pic{
\draw [fill=black] (0,0) circle [radius=.05];
\draw [fill=black] (1,0) circle [radius=.05];
\draw [fill=black] (2,0) circle [radius=.05];
\draw [fill=black] (3,0) circle [radius=.05];
\draw [fill=black] (4,0) circle [radius=.05];
\draw [fill=black] (5,0) circle [radius=.05];
\draw [fill=black] (6,0) circle [radius=.05];
\draw [fill=black] (7,0) circle [radius=.05];
\draw [fill=black] (8,0) circle [radius=.05];
\draw [fill=black] (9,0) circle [radius=.05];
\draw [fill=black] (10,0) circle [radius=.05];
\draw [fill=black] (11,0) circle [radius=.05];
\node at (0,-.5) {$x_1$};
\node at (11,-.5) {$x_s$};
\node at (6,-.5) {$x_j$};
\node at (7,-.5) {$x_{j+1}$};
\draw (0,0) -- (11,0);
\draw [fill=black] (0,-1.5) circle [radius=.05];
\draw [fill=black] (1,-1.5) circle [radius=.05];
\draw [fill=black] (2,-1.5) circle [radius=.05];
\draw [fill=black] (3,-1.5) circle [radius=.05];
\draw [fill=black] (4,-1.5) circle [radius=.05];
\draw [fill=black] (5,-1.5) circle [radius=.05];
\draw [fill=black] (6,-1.5) circle [radius=.05];
\draw [fill=black] (7,-1.5) circle [radius=.05];
\draw [fill=black] (8,-1.5) circle [radius=.05];
\draw [fill=black] (9,-1.5) circle [radius=.05];
\draw [fill=black] (10,-1.5) circle [radius=.05];
\draw [fill=black] (11,-1.5) circle [radius=.05];
\node at (0,-2) {$x_1$};
\node at (11,-2) {$x_s$};
\node at (6,-2) {$x_j$};
\node at (7,-2) {$x_{j+1}$};
\draw (0,-1.5) -- (6,-1.5);
\draw (7,-1.5) -- (11,-1.5);
\draw (11,-1.5) to [out=135,in=45] (6,-1.5) ;
}
\end{center}
Pos\'a then argues that the set $X$ of end-points created by a sequence of rotations has less than $2|X|$ neighbors. Then w.h.p. every set with fewer than $2|X|$ neighbors has size $\Omega(n)$ and from there he argued that $G_{n,Kn\log n}$ is Hamiltonian w.h.p. Several researchers realised that Pos\'a's argument could be tightened. Koml\'os and Szemer\'edi \cite{KS2} proved that if $m=n(\log n+\log\log n+c_n)/2$ then
\beq{1}{
\lim_{n\to\infty}\Pr(G_{n,m}\text{ is Hamiltonian})=\begin{cases}0&c_n\to-\infty.\\e^{-e^{-c}}&c_n\to c.\\1&c_n\to\infty.\end{cases}
}
Korsunov \cite{Kor} proved this for the case $c_n\to\infty$. Bollob\'as \cite{Boll2} proved the somewhat stronger hitting time result, as did Ajtai, Koml\'os and Szemer\'edi \cite{AKS}. By this we mean that  w.h.p. $m_2=m_\cH$ where $m_k=\min\set{m:\de(G_M)\geq k}$ and for a graph property $\cP$,\\ $m_\cP=\min\set{m:G_m\text{ has property }\cP}$ and $\cH=\set{\text{Hamiltonicity}}$. 

Equation \eqref{1} can be expressed as that for all $m$,
$$\lim_{n\to\infty}\Pr(G_{n,m}\text{ is Hamiltonian})\approx \lim_{n\to\infty}\Pr(\de(G_{n,m})\geq 2).$$
Alon and Krivelevich \cite{AK} proved a sort of converse, i.e.
$$\frac{\Pr(G_{n,m}\text{ is not Hamiltonian})}{Pr(\de(G_{n,m})< 2)}=1-o(1).$$
\subsection{Counting and packing}
With the existence question out of the way, other questions arise. The first concerns the number of distinct Hamilton cycles. Consider first the case of edge-disjoint Hamilton cycles. Bollob\'as and Frieze \cite{BF} proved the folowing. Let property $\cA_k$ be the existence of $\rdown{k/2}$ edge disjoint Hamilton cycles plus a disjoint matching of size $\rdown{n/2}$ if $k$ is odd. They proved that if $k=O(1)$ then
\beq{Ak}{
m_{\cA_k}=m_{k}\text{ w.h.p. }
}
It took some time to solve the question of dealing with the case of growing $k$. It is marginally weaker to say that w.h.p. $G_{n,p}$ has property $\cA_{\de}$ where here $\de=\de(G_{n,p})$. Frieze and Krivelevich \cite{FrKr} proved this is true as long as $np=(1+o(1))\log n$ and Ben-Shimon, Krivelevich and Sudakov \cite{BKS} extended the range to $np\leq (1.02)\log n$. And then the dam broke and Krivelevich and Samotij \cite{KSa} proved that w.h.p. $G_{n,p}$ has property $H_\de$ for $p=O(n^{-\a})$ where $\a<1$ is a positive constant and Knox, while K\"uhn and Osthus \cite{KKO} proved that $G_{n,p}$ has property $H_\de$ w.h.p. for $\log^{50}n/n\leq p\leq 1-n^{-1/4}\log^9n$.
\begin{Problem}
Is it true that w.h.p. $m_{\cA_k}=m_{k}$ holds true throughout the whole of the graph process?
\end{Problem}
Briggs, Frieze, Krivelevich, Loh and Sudakov \cite{BFKLS} showed that the $k$ disjoint Hamilton cycles can be found on-line. Let $\t_{2k}$ be the hitting time for minimum degree at least $2k$. In \cite{BFKLS} it is shown that w.h.p. the first $\t_{2k}$ edges can be partitioned on-line into $k$ subsets, so that each subset contains a Hamilton cycke.

The above results concern packing Hamilton cycles. In the dual problem, we wish to cover all the edges by a small collection of hamilton cycles. A trivial lower bound for the number of cycles needed to cover the edges of a graph $G$ is $\rdup{\D(G)/2}$, where $\D$ denotes maximum degree. Glebov, Krivelevich and Szab\'o \cite{GKS} studied expander graphs and proved that w.h.p. $(1+o(1))\D/2$ are sufficient for $G_{n,p},\,p\geq n^{-1+\e}$. Hefetz, K\"uhn, Lapinskas and Osth\"{u}s \cite{HKLO} proved the tight result, i.e.  $\rdup{\D(G)/2}$ are sufficient, for $\frac{\log^{117}n}{n}\leq p\leq  1-n^{-1/8}$. Dragani\'c, Glock, Corriea and Sudakov \cite{DGCS} extended this result to all $p\geq C\log n/n$ for large enough $C$.

We next consider the question of the number of distinct Hamilton cycles in a random graph. Let $X_H=X_H(G)$ denote the number of Hamilton cycles in the graph $G$. Janson \cite{Ja1} proved that if $m\gg n^{3/2}$ and $N-m\gg n$ then $(X_H-\E(X_H))/\Var(X_H)^{1/2}$ converges in distribution to the standard normal distribution. He also proved that $G_{n,p}$ behaves differently, in the sense that the number of Hamilton cycles converges in distribution to a log-normal distribution when $np\gg n^{1/2}$, but $p<\a<1$ for some constant $\a>0$. Normality for $G_{n,p}$ only happens for $p\to 1$.

There is still the question of how large is $X_H$ at the hitting time $m_\cH$for Hamilton cycles. Cooper and Frieze \cite{CF1} showed that w.h.p. $G_{m_\cH}$ contains $(\log n)^{n-o(n)}$ Hamilton cycles, which is best possible up to the value of the $o(n)$ term. Glebov and Krivelevich \cite{GK} proved that $(\log n)^{n-o(n)}$ can be improved to $(\log n/e)^n(1-o(1))^n$ . On the other hand, if we want the expected number of Hamilton cycles at time $m_{\cH}$ then McDiarmid \cite{McD1} proved that $\E(X_H)\approx 8(n-1)!(\p n)^{1/2}4^{-n}$. The discrepancy between this and previous results stems from the fact that the expectation is dominated by the likely number of Hamilton cycles when the hitting time is $\Omega(n^2)$. This number compensates for the unlikely hitting time of $\Omega(n^2)$.
\begin{Problem}\label{P3}
W.h.p., at time $m_{\cH}$, there are $n!p^ne^{o(n)}$ Hamitlon cycles. Determine $o(n)$ as accurately as possible. 
\end{Problem}
Glock and Sgueglia \cite{GlSg} proved that for every fixed $k$ there exists $C=C(k)$ such that if $p\geq C\log n/n$ then w.h.p. $G_{n,p}$ contains $k$-edge disjoint perfect matchings with the property that every pair of such matchings form a Hamilton cycle. 
\subsection{Lower bounds on the minimum degree}
We have seen that the threshold for Hamilton cycles is intimately connected to the threshold for minimum degree at least two. More generally, the threshold for the property $\cA_k$ is connected to the threshold for minimum degree at least $k$. So, if we condition our graphs to have minimum degree $k$ then we should have a lower threshold. Bollob\'as, Fenner and Frieze \cite{BFF1} considered the random graph $G_{n,m}^{(k)}$. This being a random graph selected uniformly from the set $\cG_{n,m}^{(k)}$ of graphs with vertex set $[n]$, $m$ edges and minimum degree at least $k$. They proved that if $m=\frac{n}{2(k+1)}\brac{\log n+k(k+1)\log\log n+c_n}$ then
\beq{Hamk}{
\lim_{n\to\infty}\Pr(G_{n,m}^{(k)}\in\cA_k)=\begin{cases}0&c_n\to-\infty\text{ slowly}.\\e^{-\th_k}&c_n\to c.\\1&c_n\to\infty.\end{cases}
}
Here $\th_k=\frac{e^{-c}}{(k+1)!((k-1)!)^{k+1}(k+1)^{k(k+1)}}$. The main obstruction to $\cA_k$ is the existence of $k+1$ vertices of degree $k$, sharing a common neighbor. Also, the restriction $c_n\to-\infty$ slowly in \eqref{Hamk} is a limitation of the model being used in that paper. It can be (almost) eliminated by a better choice of model as used in the following papers. The main obstruction to being Hamiltonian for random graphs is either having minimum degree at most one or having two many vertices of degree two. When we condition on having minimum degree three, there is no natural obstruction. Bollob\'as, Cooper, Fenner and Frieze \cite{BCFF} showed that w.h.p. the random graph $G^{(k)}_{n,c_kn},c_k=(k+1)^3,k\geq 3$ has property $\cA_k$. In particular, $G^{3}_{n,64n}$ is Hamiltonian w.h.p. The value of 64 was recently reduced to 10 in Frieze \cite{F1} and then to 2.66... by Anastos and Frieze \cite{AF19}.
\begin{Problem}\label{P4}
Is $G^{(3)}_{n,cn},c>3/2$ Hamiltonian w.h.p.? More generally, does $G^{(k)}_{n,d_k/n},d_k>k/2$ have property $cA_k$?
\end{Problem}
The paper \cite{KLS} by Krivelevich, Lubetzky and Sudakov proves that in the random graph process, the $k$-core, $k\geq 15$ has Property $\cA_{k-1}$ w.h.p., as soon as it is non-empty. Thus we immediately get the problem:
\begin{Problem}
Replace $k\geq 15$ by $k\geq 3$ and $\cA_{k-1}$ by $\cA_k$ in the result of \cite{KLS}.
\end{Problem} 
Anastos \cite{Ana22c} replaced $k\geq 15$ by $k\geq 4$.
\subsection{Resilience}
Sudakov and Vu \cite{SV} intoduced the notion of {\em (local) resilience}. In our context, the local resilience of the Hamiltonicity property is the maximum value $\D_{ham}$ so that w.h.p. $G_{n,p}-H$ is Hamiltonian for all $H\subseteq G$ with maximum degree $\D(H)\leq \D_{ham}$. The aim now is to prove a result with $\D_{ham}$ as large as possible and $p$ as small as possible. We let $\cL(p,\D)$ denote that $G_{n,p}$ has local resilience of hamiltonicity for $\D_{ham}\leq\D$. Sudakov and Vu proved local resilience for $p\geq \frac{\log^4n}{n}$ and $\D_{ham}=\frac{(1-o(1))np}{2}$. The expression for $\D_{ham}$ is best posible, but the needed value for $p$ has been lowered. Frieze and Krivelevich \cite{FrKr} showed that there exist constants $K,\a$ such that $\cL\brac{\frac{K\log n}{n},\a np}$ holds w.h.p. Ben-Shimon, Krivelevich and Sudakov \cite{BKS} improved this to $\cL\brac{\frac{K\log n}{n},\frac{(1-\e)np}{6}}$ holds w.h.p. and then in \cite{BKS2} they obtained a result on resilience for $np-(\log n+\log\log n)\to \infty$, but with $K$ close to $\frac{1}{3}$. (Vertices of degree less than $\frac{np}{100}$ can lose all but two incident edges.) Lee and Sudakov \cite{LS} proved the sought after result that for every positive $\e$ there exists $C=C(\e)$ such that w.h.p. $\cL\brac{\frac{C\log n}{n},\frac{(1-\e)np}{2}}$ holds. Condon, Espuny D\'iaz, Kim, K\"uhn and Osthus \cite{CDKKO1} refined \cite{LS}. Let $H$ be a graph with degree sequence $d_1\geq d_2\geq \cdots \geq d_n$ where $d_i\leq (n-i)p-\e np$ for $i<n/2$. They say that $G$ is $\e$-P\'osa-resilient if $G-H$ is Hamiltonian for all such $H$. Given $\e>0$ there is a constant $C=C(\e)$ such that if $p\geq \frac{C\log n}{n}$ then $G_{n,p}$  is $\e$-P\'osa-resilient w.h.p. 

The result in \cite{LS} has now been improved to give a hitting time result, see Montgomery \cite{M1} and Nenadov, Steger and Truji\'c \cite{NST}. The latter paper also proves the optimal resilience of the 2-core when $p=\frac{(1+\e)\log n}{3n}$. It would seem that the Hamiltonicity resilience problem is completely resolved, but one can still ask the following:
\begin{Problem}
Assuming that $G^{(3)}_{n,cn}$ is Hamiltonian w.h.p., what can one say about its resilience? 
\end{Problem}
 Alon \cite{Yalon} studied global resilience and showed that if $p$ is above the threshold for Hamiltonicity in $G=G_{n,p}$ then one has to remove at least $\delta(G)-1$ edges to destroy Hamiltonicity.
\begin{Problem}
How many edges must be deleted so that the fewest number of vertex disjoint paths covering $[n]$ is equal to $k$. Is it $k(\delta(G)-1)$?
\end{Problem}
\subsection{Powers of Hamilton cycles}
The $k$th power of a Hamilton cycle in a graph $G=(V,E)$ is a permutation $x_1,x_2,\ldots,x_n$ of the vertices $V$ such that $\set{x_i,x_{i+j}}$ is an edge of $G$ for all $i\in [n],j\in[k]$. K\"uhn and Osthus \cite{KO} studied the existence of $k$th powers in $G_{n,p}$. They showed that for $k\geq 3$ one could use Riordan's Theorem \cite{R} to show that if $np^k\to\infty$ then $G_{n,p}$ contains the $k$th power of a Hamilton cycle w.h.p. For $k=2$ they only showed that $np^{2+\e}\to\infty$ was sufficient. Subsequently Nenadov and \v{S}kori\'c \cite{NS} showed that if $np^2\geq C\log^8n$ for sufficiently large $C$ then $G_{n,p}$ contains the square ($k=2$) of a Hamilton cycle w.h.p. Fischer, \v{S}kori\'c, Steger and Truji\'c \cite{FSST} have shown that there exists $C>0$ such that if $p\geq \frac{C\log^3n}{n^{1/2}}$ then not only is there the square of a Hamilton cycle w.h.p., but containing a square is resilient to the deletion of not too many triangles incident with each vertex. Montgomery \cite{Mnotes} improved the bound to $p\gg \frac{\log^2n}{n^{1/2}}$. 

It is interesting that P\'osa rotations have played a significant role in everything mentioned so far, except for \cite{NS}. They used the {\em absorbing} method, and this plays a role in other recent papers. As discussed in \cite{NS}, we can demonstrate the basic idea in the simpler case of Hamilton cycles in graphs. Let $A$ be a graph and $a,b\in V(A)$ two distinct vertices. Given a subset $X\subseteq V(A)$, we say that $A$ is an $(a,b,X)$-absorber if for every subset $X'\subseteq X$ there exists a path $P_{X'}\subseteq A$ from $a$ to $b$ such $V(P)=V(A)\setminus X'$. Let $G=(V,E)$ be a graph in which we want to find a Hamilton cycle and suppose there exists a large subset $X\subseteq V$ and an $(a,b,X)$-absorber $A\subseteq G$, for some vertices $a,b\in V(A)$. An important observation is that if $G$ contains a path from $a$ to $b$ such that
$P$ uses all the vertices in $V\setminus V(A)$ and no vertex from $V(A)\setminus X$ (except $\set{a,b}$), we are done. Indeed, if $X'$ is the subset of $X$ used by $P$ then by the definition of absorber, there is a path $P_{X'}\in A$ which together with $P$ gives a Hamilton cycle. 

It takes work to show the existence of $P$ and absorbers, but it is definitely introduces a new idea to Hamilton cycle problems in random structures. 

More recently, a general result on thresholds by Frankston, Kahn, Narayanan and Park \cite{FKNP} improves the bound to $p\geq K\frac{\log n}{n^{1/2}}$ for sufficiently large $K$. A refinement of this approach in  Kahn, Narayanan and Park \cite{KNP} improves the bound to $p\geq \frac{K}{n^{1/2}}$ for sufficiently large $K$. It was conjectured in \cite{KNP} that $K\approx e^{-1/2}$. Finally, Zhukovskii \cite{Zuk} showed that the threshold for a $k$th power is asymptotic to $(e/n)^{1/k}$,  for $k\geq 2$, resolving the conjecture in \cite{KNP}. Makai, Pasch, Petrova, Schiller \cite{MPPS} independently verified the same threshold but only for $k\geq 4$.  
\subsection{Edge-colored Random Graphs}
 Many nice problems arise from considering random graphs with colored edges.
 \subsubsection{Rainbow Hamilton Cycles}
 A set of colored edges $E$ is called {\em rainbow} if every edge has a different color. Cooper and Frieze \cite{CF2} proved that if $m\geq 21n\log n$ and each edge of $G_{n,m}$ is randomly given one of at least $q\geq 21n$ random colors then w.h.p. there is a rainbow Hamilton cycle. Frieze and Loh \cite{FLo} improved this result to show that if $m\geq \frac{1}{2}(n+o(n))\log n$ and $q\geq (1+o(1))n$ then w.h.p. there is a rainbow Hamilton cycle. This was further improved by Ferber and Krivelevich \cite{FeKr} to $m=n(\log n+\log\log n+\om)/2$ and $q\geq (1+o(1))n$, where $\om\to\infty$ with $n$. This is best possible in terms of the number of edges.
 \begin{Problem}\label{P9}
 Suppose that $q=cn,c\geq 1$ and that we consider the graph process $G_0,G_1,\ldots,G_m$. Let
\beq{hit}{
   \t_c=\min\set{t:G_t\text{ contains $n$ distinct edge colors}}\text{ and }\t_2=\min\set{t:\delta(G_t)\geq 2}.
}
Is it true that w.h.p. there is a rainbow Hamilton cycle at time $\max\set{\t_c,\t_2}$?\\
Frieze and McKay \cite{FM} proved the equivalent of \eqref{hit} when Hamilton cycle is replaced by spanning tree. (Here we required $\delta(G_t)\geq 1$.)
\end{Problem}
The case $q=n$ was considered by Bal and Frieze \cite{BalF}. They showed that $O(n\log n)$ random edges suffice.
\begin{Problem}\label{P10}
Discuss the problem of packing rainbow Hamilton cycles in $G_{n,m}$. Are there rainbow colored versions of \cite{BF}, \cite{KSa} and \cite{KO}? Ferber and Krivelevich \cite{FeKr} give asymptotic results along this line.
\end{Problem} 
\begin{Problem}
Typically, how many rainbow Hamilton cycles are there in $G_{n,p}$ for a given number of colors $q$?
\end{Problem}
Ferber, Han and Mao \cite{FHM} used a different definition of rainbow cycles. Here $G_1,G_2,\ldots,G_n$ are independent copies of $G_{n,p}$. A rainbow Hamilton cycle is one made up of $n$ edges, one from each $G_i$. They prove a Dirac type result. They show that if $H_i\subseteq G_i$ and $\de(H_i)\geq (\frac12+\e)np$ for $i=1,2,\ldots,n$ and that if $p\gg\bfrac{\log n}{n}$ then w.h.p. the family $H_1,H_2,\ldots,H_n$ admits a rainbow Hamilton cycle.

Anastos and Chakraborti \cite{AnCh} considered two similar scenarios. In the first we have a set of graphs $G_i\cap G_{n,p}$ where $G_i$ is Dirac graph, i.e. $\delta(G_i)\geq n/2$, for $i=1,2,\ldots,n$ and $p\geq c\log n/n$. In the second, $G$ is a fixed Dirac graph and $G_i$ is an independent copy of $G\cap G_{n,p},p\geq c\log n/n^2$. In both cases there is a rainbow hamilton cycle w.h.p.

Christoph, Martin and Milojevic \cite{CMM2025} proved that if $p\geq C\log n/n$ and $\cG=(G_1,G_2,\ldots,C_n)$ are independent copies of $G_{n,p}$ then w.h.p. $\cG$ is Hamilton-Universal. By this we mean that for every map $f:[n]\to[n]$ there is a Hamilton cycle whose $i$th edge comes from $G_{f(i)}$.
\begin{Problem}
Find the optimal constant for Hamilton-Universality.
\end{Problem}
\subsection{Anti-Ramsey property}
The rainbow concept is closely related to the Anti-Ramsey concept. Introduced by Erd\H{o}s, Simonovits and S\'os \cite{ESS}. Cooper and Frieze \cite{CF3} considered the following. Suppose we are allowed to color the edges of $G_{n,p}$, but we can only use any color $k=O(1)$ times, a {\em $k$-bounded coloring}. They determined the threshold for every $k$-bounded coloring of $G_{n,p}$ to have a rainbow Hamilton cycle.
\begin{Problem}
Remove the upper of $O(1)$ in \cite{CF3}. Consider the case where the bound only applies to the edges incident with the same vertex. Consider the case wher the coloring is proper.
\end{Problem}
\subsection{Pattern Colorings}\label{PC}
Given a coloring of the edges of a graph, there are other patterns that one can search for with respect to Hamilton cycles. For example Espig, Frieze and Krivelevich \cite{EFK} considered {\em zebraic} Hamilton cycles. Here the edges of $G_{n,p}$ are randomly colored black and white. A Hamilton cycle is  zebraic if its edges alternate in color. They showed that the hitting time for the existence of a zebraic Hamilton cycle coincides with the hitting time for every vertex to be incident with an edge of both colors. They related this to the question of how many random edges must be added to a fixed perfect matching $M$ of $K_n$ so that there exists a Hamilton cycle $H$ that contains $M$. This turns out to coincide with the number of edges needed for minimum degree one.

Suppose next that we have used $r$ colors to randomly color edges and we have a fixed pattern $\Pi$ of length $\ell$ in mind. We say that a Hamilton cycle with edges $e_1,e_2,...,e_n$ is $\Pi$-colored if $e_j$ has color $\Pi_{t}$, where $t=j\mod\ell$. It is shown by Anastos and Frieze \cite{AF} that w.h.p. the hitting time for the existence of a $\Pi$-colored Hamilton cycle coincides with the hitting time for every vertex to {\em fit} $\Pi$. We say that vertex $v$ fits $\Pi$ if there exists $1\leq j\leq \ell$ and edges $f_1,f_2$ incident with $v$ such that $f_1$ has color $\Pi_{j}$ and $f_2$ has color $\Pi_{j+1}$.
\begin{Problem}\label{PC1}
Find precise thresholds for rainbow colored powers of Hamilton cycles or pattern colored Hamilton cycles.
\end{Problem}
Bell and Frieze \cite{BellFrieze} established these thresholds up to a constant factor.

In problems \ref{P9} and \ref{PC1} we have assumed that colors are chosen uniformly.
\begin{Problem}
Modify problems \ref{P9} and \ref{PC1} by assuming that color $c$ is chosen with probability $p_c$, for $c\in C$, the set of available colors.
\end{Problem}
Frieze and Pegden \cite{arbpattern} consider more general patterns. A color pattern will be a sequence $\bc=(c_1,c_2,\ldots,c_n)$. Given a sequence \bc\ we say that the Hamilton cycle $H=(x_1,x_2,\ldots,x_n,x_1)$ (as a sequence of vertices) is \bc-colored if $c(\set{x_i,x_{i+1}})=c_i$ for $i=1,2,\ldots,n$. 

Suppose that $\va=(\a_1,\a_2,\ldots,\a_k)$ where $\a_1,\ldots,\a_k$ are constants and $\a_1+\cdots+\a_k=1$ and $\a_i>0,i=1,2,\ldots,k$. Let $\b=\min\set{\a_i:i\in[k]}^{-1}$ and let $G_{n,p;\va}$ denote the random graph $G_{n,p}$ where each edge is independently given a random color $i$ from the {\em palette} $[k]$ with probability $\a_i$. They show that w.h.p. if $p=(\log n+\om)/n$ where $\om\to\infty$ then w.h.p. $G_{n,\b p;\va}$ contains a \bc-colored Hamilton cycle.
\begin{Problem}
Is the same probability $p$ sufficient to imply that  $G_{n,\b p;\va}$  {\em simultaneously} contains \bc-colored Hamilton cycles for all possibe \bc.
\end{Problem}
The paper \cite{arbpattern} also consides the case where the vertices are colored. Replacing edges by vertices in \bc-colored they show that $p=K\log n/n$ is sufficient for the existence of a \bc-colored Hamilton cycle w.h.p.
\begin{Problem}
Can $K$ be replaced by $\b$ as in the edge case?
\end{Problem}
\subsection{A Ramsey type question}
Gishboliner, Krivelevich and Michaeli \cite{GKM} considered the following problem. Given an $r$ coloring of $G_{n,p}$ what is the maximum number of same colored edges can we find in a Hamilton cycle. They prove that with $p$ above the Hamiltonicity threshold, in any $r$-colouring of the edges there exists a Hamilton cycle with at least $((2-o(1))/(r+1))n$ edges of the same colour. This estimate is asymptotically optimal.
\begin{Problem}
What can be said about the colorings of edge disjoint Hamilton cycles in this framework?
\end{Problem}
\subsection{Color profile}
Chakraborty, Frieze and Hasabanis \cite{CFH} considered the following problem. For a fixed $r\geq 1$ let $\bM$ denote the set $\set{{\bm}=(m_1,m_2,\ldots,m_r)\in [0,n]^r:m_1+\cdots+m_r=n}$. For a graph $G$ with $n$ vertices and edges colored with $r$ colors, let $hcp(G)$ denote the set of $\bm\in\bM$ such that $G$ contains a Hamilton cycle that is the concatenation of paths $P_1,P_2,\ldots,P_r$ such that $P_i$ contains $m_i$ edges colored $i$, for $i=1,2,\ldots,r$. The paper \cite{CFH} considers $hcp(G_{n,p})$ when the coloring is random and either (i) $p=\frac{\log n+r\log\log n+\om}{n}$ or (ii) $p=\frac{\log n+\log\log n+\om}{\a_{\min}n}$. Here we use color $i$ with probability $\a_i$. In case (i) they get a good estimate of $hcp(G_{n,p})$ close to the Hamiltonicity threshold and in (ii) they show that $hcp(G_{n,p})=\bM$ w.h.p.
\begin{Problem}
(a) Determine $hcp(G_{n,p})$ w.h.p. for all $p\leq \frac{\log n+\log\log n+\om}{\a_{\min}n}$. \\
(b) Determine $hcp(G)$ for other models of a random graph e.g. random regular graphs, or random geometric graphs.\\
(c) Extend the notion to random hypergraphs.
\end{Problem}
\subsection{Perturbations of dense graphs}\label{pert}
Spielman and Teng \cite{ST} introduced the notion of {\em smoothed analysis} in the context of Linear Programming. This inspires the following sort of question. Suppose that $H$ is an {\em arbitrary} graph and we add some random edges $X$, when can we assert that the graph $G=H+X$ has some particular property? The first paper to tackle this question was by Bohman, Frieze and Martin \cite{BFM} in the context of Hamiltonicity. They show that if $H$ has $n$ vertices and its minumum degree is at least $dn$ for some positive constant $d\leq 1/2$ and $|X|\geq 100n\log d^{-1}$ then $G$ is Hamiltonian w.h.p. This is best possible in the sense that there are bipartite graphs with minimum degree $dn$ such that adding less than $\frac13 n\log d^{-1}$ edges leaves a non-Hamiltonian graph w.h.p. Ma and Yan \cite{MY} replace 100 by $(1+\e)$ and show that this is best possible.

Further, with an upper bound on the size of an independent set, we only need $|X|\to\infty$ when $d$ is constant. Espuny D\'iaz and Razafindravola \cite{DR} considered the case where $d=n/2-\eta$ and show that if $\eta=o(n)$ and $\eta=\om(1)$ then $(8+o(1))\eta$ random edges suffice w.h.p.

 Dudek, Reiher, Ruci\'nski and Schacht \cite{DRRS} proved that if the minimum degree of  $H$ is at least $\a>k/(k+1$ then w.h.p. $H$ plus $O(n)$ random edges yields a graph containing the $(k+1)$th power of a Hamilton cycle. Nenadov and Truji\'c \cite{NeTr} improved this by showing that under the same conditions there is also a $(2k+1)th$ power. Antoniuk, Dudek and Ruci\'nski \cite{ADS} considered the case where $H$ has minimum degree at least $(1/2+\e)n$. They prove further results in \cite{ADS1}.
\begin{Problem}  
Can the constructions in \cite{DRRS} or \cite{ADS} be done in polynomial time? 
\end{Problem}B\"ottcher, Montgomery, Parczyk and Person \cite{BMPP} show that for each $k\geq 2$ there is some $\eta=\eta(k,\a)>0$ such that if $H$ has minimum degree at least $\a n$ and $|X|\geq n^{2-1/k-\eta}$ then w.h.p. $H+X$ contains a copy of the $k$th power of a Hamilton cycle. Note that $m=n^{2-1/k}$ is the threshold number of edges that $G_{n,m}$ needs for the $k$th power of a Hamilton cycle. At least for $k\geq 3$. For $k=2$ there is still a $\log n^{O(1)}$ factor to be removed. Antoniuk, Dudek, Reiher, Ruci\'nski and Schacht \cite{ADRRS} show that that adding $O(n^{2-2/\ell})$ random edges to an $n$-vertex graph G with minimum degree at least $\a n$ yields the existence of the $(k\ell+r)$-th power of a Hamiltonian cycle w.h.p.
\begin{Problem}
Determine the best possible value of $\eta$ in \cite{BMPP}.
\end{Problem}
Espuny Diaz and Razafindravola \cite{DR} consider adding random edges to a graph $H$ with minimum degree $\delta(H)$. Suppose that $\eta=n/2-\delta(H)$. They show that if $\eta=\omega(1)$ then $O(\eta)$ edges suffice to make the graph Hamiltonian w.h.p. If $\eta=O(1)$ then $\omega(1)$ edges suffice. 

Anastos and Frieze \cite{AF1} considered the addition of $m$ randomly colored edges $X$ to a randomly edge colored dense graph $H$ with with minimum degree at least $\delta n$.  The colors are chosen randomly from $[r]$ and $\th=-\log\delta$. They show that if $m\geq \min \set{ (435+75\theta)tn,  \card{\binom{[n]}{2} \setminus E(H)} }$ and $r\geq(120+20\theta)n$ then,  w.h.p. $H+X$ contains $t$ edge disjoint rainbow Hamilton cycles. Aigner-Horev and Hefetz \cite{AigH} showed that adding $r=n+o(n)$ random colors are suffiicient for Rainbow Hamitonicity. Kasamaktsis and Letzter \cite{KaLe} proved the same result for the case of $r=n$ colors.
\begin{Problem}
Determine the exact relationship between the various parameters $\de,C,r$ that ensure the existence of a rainbow Hamilton cycle w.h.p.
\end{Problem}
Espuny D\'iaz and Gir$\tilde{\text{a}}$o \cite{DG} consider the effect of adding a random $r$-regular graph $H$ to a dense graph $G$, $r\in\set{1,2}$. When $r=2$ they show that $G\cup H$ is pancyclic w.h.p. for any $G$ with minimum degree $dn, d>0$. When $r=1$ they show that $d\in [\sqrt{2}-1]$ is necessary and sufficient. Dragani\'c and Keevash \cite{DK} found a tight threshold of $\sim (n\log n/2)^{1/2}$ for the case $r=2$. Henderson, Longbrake, Mao and Morawski \cite{HLMM} proved that for each $d\geq 2$ the union of an arbitrary  $d$-regular graph and a random 2-factor is Hamiltonian w.h.p. 

Espuny D\'iaz \cite{ED} studied the effect of adding a random geometric graph to a dense graph. He showed that w.h.p. adding random geometric graph with radius $(C/n)^{1/d}$ to a graph with minimum degree at least $\a n$ is sufficient to get a hamilton cycle w.h.p. Espuny D\'iaz and Hyde \cite{EH} consider powers of Hamilton cycles in this model.

Katsamaktis, Lezter and Sgueglia \cite{KaLS} prove a general theorem relating the existence of rainbow copies in randomly colored perturbations of dense graphs to the existence of uncolored copies in randomly perturbed dense graphs.
\subsection{Compatible Cycles}
Given a graph $G=(V,E)$, a {\em compatibility system} is a family $\cF=\set{F_v:v\in V}$ of sets of edges. Each $F_v$ consists only of edges incident with $v$. The incompatibility system is $\m n$ bounded if $|F_v|\leq \m n$ for all $v\in V$. A Hamilton cycle is compatible with $\cF$ if it uses at most one edge from each $F_v,v\in V$. Krivelevich, Lee and Sudakov \cite{KLS1} proved that there exists $\m>0$ such that if $p\gg \frac{\log n}{n}$ then w.h.p. $G_{n,p}$ contains a compatible Hamilton cycle for every $\m n$ bounded compatibility system.
\begin{Problem}
Determine the maximum value of $\m>0$ for which $G_{n,p},p\gg \frac{\log n}{n}$ contains a compatible Hamilton cycle for every $\m n$ bounded compatibility system w.h.p. (The bound in \cite{KLS1} is small, but increases to $1-\frac{1}{\sqrt{2}}$ for $p\gg \frac{\log^8n}{n}$.
\end{Problem}
\subsection{Cycle space}
Christoph, Nenadov and Petrova \cite{CNP} show that if $np\geq C\log n$ for large enough $C$ then w.h.p. the Hamilton cycles of $G_{n,p}$ span its cycle space. Hefetz and Krivelevich \cite{HefKriv} improved this to a best possible $np=\log n+2\log\log n+\om(1)$. Hefetz and Krivelevich \cite{HefKriv1} extended this property to random $d$-regular graphs $G_{n,d}, n$ odd and $d$ sufficiently large and also to the randomly perturbed model of \cite{BFM}.
\subsection{Algorithms}
Finding a Hamilton cycle in a graph is an NP-hard problem. On average, however, things are not so bleak. Angluin and Valiant \cite{AV} gave a polynomial time randomised algorithm that finds a Hamilton cycle w.h.p. in $G_{n,p}$ for $p\geq \frac{K\log n}{n}$ when $K$ is sufficiently large. Shamir \cite{S} gave a  polynomial time randomised algorithm that finds a Hamilton cycle w.h.p. if $p\geq \frac{\log n+(3+\e)\log\log n}{n}$. Bollob\'as, Fenner and Frieze \cite{BFF} gave a deterministic $O(n^{3+o(1)})$ time algorithm HAM with the property 
\[
\lim_{n\to \infty}\Pr(HAM\text{ finds a Hamilton cycle in  }G_{n,m})=\lim_{n\to \infty}\Pr(G_{n,m} \text{contains a Hamilton cycle.})
\]
Nenadov, Steger and Su \cite{NSU} gave an $O(n)$ time randomised algorithm that succeeds w.h.p. when $m\geq Cn\log n$ for $C$ sufficently large. They pose the following question:
Is there an $O(n)$ time algorithm that succeeds w.h.p. at the threshold for Hamiltonicity? This was answered in the affirmative by Anastos \cite{Ana22b}.

The above algorithms used extensions and rotations. For dense random graphs Gurevich and Shela \cite{GS} gave a simpler randomised algorithm that determines the Hamiltonicity of $G_{n,p}$ in $O(n^2)$ expected time for $p$ constant. Here the algorithm resorts to using the Dynamic Programming algorithm of Held and Karp \cite{HK} if it fails to find a Hamilton cycle quickly in $G_{n,1/2}$. This results was strengthened to work in $G_{n,p},p\geq Kn^{-1/3}$ by Thomason \cite{T} and more reccently to $p\geq 70n^{-1/2}$ by Alon and Krivelevich \cite{AK1}.
\begin{Problem}
Can the Hamiltonicity of $G_{n,m}$ be determined in polynomial expected time for all $0\leq m\leq \binom{n}{2}$?
\end{Problem}
Anastos \cite{Ana22a} gives a deterministic algorithm for solving the Hamilton cycle problem in $G_{n,m}$. It runs in polynomial expected time provided $m\geq Kn$ for $K$ sufficiently large. Of course for $m=O(n)$ it provides a certificate of non-Hamiltonicity. 

Frieze and Haber \cite{FH} studied the algorithmic question in relation to $G_{n,cn}^{(3)}$ and showed that w.h.p. a Hamilton cycle can be found in $O(n^{1+o(1)})$ time if $c$ is a sufficiently large constant.
\begin{Problem}
Is there an $O(n\log n)$ time algorithm that w.h.p. finds a Hamilton cycle in $G_{n,cn}^{(3)}$.
\end{Problem}
One can also attack the algorithmic problem from a parallel perspective. Frieze \cite{F2} gave a parallel algorithm that uses a PRAM with $O(n\log^2n)$ processors and takes $O(\log\log n)^2$ rounds w.h.p. to find a Hamilton cycle in $G_{n,p}$, $p$ constant. MacKenzie and Stout \cite{MS} reduced the number of processors needed to $n/\log^*n$ and the number of rounds to $O(\log^*n)$. 
\begin{Problem}
Is there a PRAM algorithm that uses a polynomial number of processors and polyloglog (or better) rounds and finds a Hamilton cycle in $G_{n,p}$ at the threshold for Hamiltonicity?
\end{Problem}
In the case of Distributed Algorithms, Levy, Louchard and Petit \cite{LLP} gave an algorithm that finds a Hamilton cycle w.h.p. provided $p\gg \log^{1/2}n/n^{1/4}$. This algorithm only requires $n^{3/4+\w}$ rounds. This was recently improved to $p\gg \log^{3/2}/n^{1/2}$ in $O(\log n)$ rounds by Tureau \cite{Tu}.
\begin{Problem}
Reduce the requirements on $p$ for the existence of a distributed algorithm for finding a Hamilton cycle in $G_{n,p}$ in a sub-linear number of rounds w.h.p.
\end{Problem}
Ferber, Krivelevich, Sudakov and Vieira \cite{FKSV} considered how many edge queries one needs to find a Hamilton cycle in $G_{n,p}$. They showed that if $p\geq \frac{\log n+\log\log n+\om}{n}$ then w.h.p. one only needs to query $n+o(n)$ edges. 
\subsection{$G_p$}
Given a graph $G$ and a probability $p$, the random subgraph $G_p$ is obtained by including each edge of $G$ independently with probability $p$. A Dirac graph is a graph on $n$ vertices that has minimum degree $\delta(G)\geq n/2$. Krivelevich, Lee and Sudakov \cite{KLS3} showed that if $G$ is a Dirac graph and $p\geq \frac{C\log n}{n}$ then $G_p$ is Hamiltonian w.h.p. Given a graph $G$ we can define a graph process $G_0,G_1,\ldots,$ where $G_{m+1}$ is obtained from $G_m$ by adding a random edge from $E(G)\setminus E(G_m)$. Johansson \cite{Jo} showed that if $G$ has minimum degree at least $(1/2+\e)n$ for some positive constant $\e$ then w.h.p. the hitting time for Hamiltonicity coincides with the hitting time for minimum degree at least two. Alon and Krivelevich \cite{AKnew} showed that w.h.p. the hitting time for $\cA_{2k}$ the same as the hitting time for minimum degree $2k$ for $k=O(1)$. They also showed this for two classes of pseudo-random graphs. Glebov, Naves and Sudakov \cite{GNS} proved that if $\delta(G)\geq k$ and $p\geq \frac{\log k+\log\log k+\om_k}{k}$ then w.h.p. (as $k$ grows) $G_p$ has a cycle of length $k+1$. When $G=K_n$ this gives part of equation \eqref{1}. Wang \cite{Wang} proved a bipartite version of \cite{Jo}.
\begin{Problem}\label{Johan}
In  \cite{AKnew}, is the hitting time for $\cA_k$ the same as the hitting time for minimum degree $k$ w.h.p. when $k=k(n)\to\infty$?
\end{Problem}
 Kelly, M\''uyesser and Pokrovskiy \cite{KMP} and Joos, Lang and Sanhueze-Matamala \cite{JLS} prove that if $\delta(G)\geq (\frac{\ell}{\ell+1}n+\e$ and $p\gg n^{-1/\ell}$ then w.h.p. $G_p$ contains the $\ell$th power of a Hamilton cycle.
\section{Random Regular Graphs}
Let $G_{n,r}$ denote a random simple regular graph with vertex set $[n]$ and degree $r$. Some of the results for $G_{n,m},G_{n,p}$ have been extended to this model.
\subsection{Existence}
Bollob\'as \cite{BolRegHam} and Fenner and Frieze \cite{FFregHam} used extensions and rotations to prove that w.h.p. $G_{n,r}$ is Hamiltonian for $r_0\leq r=O(1)$. The smaller value of $r_0$ here was 796. At around the same time Robinson and Wormald \cite{RW} showed that random cubic bipartite graphs are Hamiltoninan w.h.p. The gap for $3\leq r=O(1)$ was filled by Robinson and Wormald \cite{RW1}, \cite{RW1a}. They introduced an ingenious variation on the second moment method that is now referred to as small subgraph conditioning. Basically, it says, in some sense, that if we condition on the number of small odd cycles then the second moment method will prove that $G_{n,r}$ is Hamiltonian w.h.p. 

This leaves the case for $r\to\infty$. In unpublished work, Frieze \cite{F0} proved that $G_{n,r}$ is Hamiltonian w.h.p. for $3\leq r\leq n^{1/5}$. This was improved to $r\leq c_0n$ for some constant $c_0>0$ by Cooper, Frieze and Reed \cite{CFR}. At the same time Krivelevich, Sudakov, Vu and Wormald \cite{KSVW} proved the same result for $r\geq n^{1/2}\log n$.

Frieze \cite{Fperm} proved that w.h.p. the union of two random permutation graphs on $[n]$ contains a Hamilton cycle. Here we ignore orientation. This has some relation to Theorem 4.15 of \cite{Wormreg}.

Robinson and Wormald \cite{RWmat} show that we can specify $o(n^{1/2})$ edges of a matching $M$, with an orientation, and w.h.p. find a Hamilton cycle $H$ in $G_{n,r},r\geq 3$ that contains $M$. Here $H$ the edges of $M$ appear on $H$ with the correct orientation. This implies that a random {\em claw-free} cubic graph is Hamiltonian w.h.p. The same paper also shows that if $|M|=o(n^{2/5})$ then we can impose an ordering of the edges around the cycle.

Kim and Wormald \cite{KW} showed that w.h.p. $G_{n,r},r\geq 3$ satisfies property $\cA_r$. Thus w.h.p. $G_{n,r}$ is the union of edge disjoint Hamilton cycles and a perfect matching if $r$ is odd.
\subsection{Algorithms}
Frieze \cite{Freg} showed that the extension-rotation approach gives rise to an $O(n^{3+o(1)})$ time algorithm that finds a Hamilton cycle in $G_{n,r},85\leq r=O(1)$ w.h.p. Frieze, Jerrum, Molloy, Robinson and Wormald \cite{FJMRW} found an approach that works for $r\geq 3$. It follows from the work of Robinson and Wormald \cite{RW1}, \cite{RW1a} that w.h.p. the number of 2-factors of $G_{n,r}$ is at most $n$ times the number of Hamilton cycles in $G_{n,r}$. So, if we generate a near uniform 2-factor, it has probability of a least $n^{-1}$ of being a Hamilton cycle. If we generate $n\log n$ random 2-factors, then w.h.p. one of them will be a Hamilton cycle. To generate a random 2-factor, we use the Markov chain approach of Jerrum and Sinclair \cite{JS}.
\begin{Problem}\label{P21}
Construct a near linear time algorithm for finding a Hamilton cycle in $G_{n,r}, r\geq 3$.
\end{Problem}
\subsection{Rainbow Hamilton Cycles}
Janson and Wormald \cite{JW} proved the following: Suppose that the edges of of the random $2r$-regular graph $G_{n,2r}$ are randomly colored with $n$ colors so that each color is used exactly $r$ times. Then w.h.p. there is a rainbow Hamilton cycle if $r\geq 4$ and there isn't if $r\leq 3$.
\begin{Problem}
Discuss the problem of packing rainbow Hamilton cycles in the context of random regular graphs.
\end{Problem}
\subsection{Resilience}
Condon, Espuny D\'iaz, Gir$\tilde{\text{a}}$o, K\"uhn and Osthus \cite{CDGKO} proved that given $\e>0$, $\D_{ham}\leq (\frac12-\e)d$ for $d_\e\leq d\leq \log^2n$, and the lower bound is necessary. The result for larger values of $d$ follows by combining results from several papers. Sudakov and Vu \cite{SV} showed that for any fixed $\e>0$, and for any $(n,d,\l)$-graph $G$ with $d/\l>\log^2n$, we have that $\D_{ham}\leq (\frac12-\e)d$. This, together with a result of Krivelevich, Sudakov, Vu and Wormald \cite{KSVW} and recent results of Cook, Goldstein and Johnson \cite{Cook} and Tikhomirov and Youssef \cite{Tik} about the spectral gap of random regular graphs implies that $\D_{ham}\leq (\frac12-\e)d$ for $\log^4n\leq d\leq n-1$ w.h.p. One can extend this to $d\gg \log n$ by combining the coupling result of Kim and Vu \cite{KimVu} with that of Lee and Sudakov \cite{LS}.

\subsection{Fixed Degree Sequence}
Regualrity is a simple example of a fixed degree sequence. Let $\bd=(d_1,d_2,\ldots,d_n)$ be a degree sequence. We let $G_{n,\bd}$ denote a random graph chosen uniformly from all graphs with vertex set $[n]$ and with degree sequence $\bd$. Cooper, Frieze and Krivelevich \cite{CFK} gave some rather complicated conditions on \bd\ under which $G_{n,\bd}$ is Hamiltonian w.h.p. Gao, Isaev and McKay \cite{GIM} proved Hamiltonicity w.h.p. under the assumption that \bd\ is {\em near-regular} viz. $\max d_i-\min d_i=o(\max d_i,n-\max d_i)$. Johansson \cite{Jo1} related Hamiltonicity to a notion of balancedness.
\begin{Problem}
Study the Hamiltonicity of $G_{n,\bd}$. Is there some simple function $\f$ such that $G_{n,\bd}$ is Hamiltonian w.h.p. if and only if $\f(\bd)>0$? (Part of the problem is to make this statement precise.)
\end{Problem}
\subsection{Sequentially Constrained Cycles}
Robinson and Wormald \cite{RW2} considered random regular graphs and ask for Hamilton cycles that contain a prescibed number of $o(n^{2/5})$ edges must be contained in order in the cycle. Espig, Frieze and Krivelevich \cite{EFK} proved that if $M$ denotes some fixed perfect matching of $K_n$ then w.h.p. at the hitting time for Hamiltonicity, there is a Hamilton cycle containing $M$ as a subgraph. Let $s_0=\om n/\log n$ where $\om=o(\log\log\log n)$ and let $S_0$ be an arbitrary $s_0$-subset of $[n]$. Frieze and Pegden \cite{arbpattern} showed that if $p=(\log n+\log\log n+\om)/n,\,\om=o(\log\log n)$ then w.h.p. $G_{n,p}$ contains a Hamilton cycle in which the vertices $S_0$ appear in natural order.
\begin{Problem}
Can $\log\log\log n$ be replaced by $\log\log n$ in the contraint on $\om$.
\end{Problem}
For a Hamilton cycle $H$ treated as a sequence of vertices $\bi=(i_1=1,i_2,\ldots,i_n)$ let $\iota(H)=|\set{k<\ell:i_k>i_\ell}|$ be the number of transpositions. Suppose that $M=\Omega(n\log n)$. Frieze and Pegden show that there is a constant $K$ such that if $p\geq \frac{Kn\log n}{M}$ then w.h.p. $G_{n,p}$ contains a Hamilton cycle $H$ with $\iota(H)\leq M$. Furthermore, if $p\leq (1-\e)\min\set{\frac{\log n}{n},\frac{n}{eM}}$ then w.h.p. $G_{n,p}$ contains no such Hamilton cycle. Here $\e$ is an arbitrary positive constant.
\section{Other Models of Random Graphs}
\subsection{Random Bipartite and multi-partite Graphs}
In the random bipartite graph $G_{n,n,p}$ we have two disjoint sets $A,B$ of size $n$ and each of the $n^2$ possible edges is included with probability $p$. Frieze \cite{Fbip} proved that if $p=\frac{\log n+\log\log n+c_n}{n}$ then
\[
\lim_{n\to\infty}\Pr(G_{n,n,p}\text{ is Hamiltonian})=\begin{cases}0&c_n\to-\infty.\\e^{-2e^{-c}}&c_n\to c.\\1&c_n\to\infty.\end{cases}
\]
Bollob\'as and Kohayakawa \cite{BoKo} proved a hitting time version and sketched a proof of the extension to $\cA_k$.
\begin{Problem}
Discuss the number of Hamilton cycles at the hitting time for minimum degree at least two in $G_{n,n,p}$.
\end{Problem}
\begin{Problem}
Discuss the resilience of Hamiltonicity in the context of random bipartite graphs.
\end{Problem}
Anastos, Frieze and Gao \cite{AFG} consider the Stochastic Block Model where the vertex set $[n]$ is partitioned into disjoint blocks $B_1,B_2,\ldots,B_k, k=O(1)$. The edge probabilities are $p$ within blocks and $q$ between blocks. They show that under some fairly general conditons on block sizes that Hamiltonicity and minimum degree two are intimately related.
\begin{Problem}
Consider the random block model where the edge probabilities are $p_{i,j}$ betwen blocks $i,j$ for $1\leq i\leq j\leq k=O(1)$.
\end{Problem}
In the vein of this problem, Johansson \cite{J2020} considered the the case where the probability of an edge between vertices $u,v$ is given by a symmetric matrix $P$. Denote this model by $G_{n,P}$. This is the ultimate version of the block model where the blocks have size one. Usually referred to as the Inhomogeneous model. He gives conditions on $P$ and shows under these conditions that w.h.p. $G_{n,P}$ has property $\cA_k$ for $k\geq 1$. He also proves hitting time results.
\subsection{$G_{k-out}$}\label{kout}
The random graph $G_{k-out}$ is a simple model of a sparse graph that has a guarnteed minimum degree. Each vertex $v\in [n]$ independently chooses $k$ random neighbors. Fenner and Frieze \cite{FeFr1} showed that $G_{23-out}$ is Hamiltonian w.h.p. Then Frieze \cite{Freg} gave a constructive proof that $G_{10-out}$ is Hamiltonian. This was followed by Frieze and {\L}uczak \cite{FL} who showed that $G_{5-out}$ is Hamiltonian. It follows from Cooper and Frieze \cite{CoFr1} that $G_{4-out}$ is Hamiltonian and then finally Bohman and Frieze \cite{BoF} showed that $G_{3-out}$ is Hamiltonian. It is easy to see that $G_{2-out}$ is non-Hamiltonian w.h.p. There must be three vertices of degree two with a common neighbor.
\begin{Problem}
Give a constructive proof that $G_{3-out}$ is Hamiltonian w.h.p.
\end{Problem}
\begin{Problem}
Does $G_{k-out},k\geq 2$ have property $\cA_{k-1}$ w.h.p. (The answer is yes, for $k=2,3$.)
\end{Problem}
There is a refinement of $G_{k-out}$ that we believe is interesting. We will call it $H_{k-out}$ where $H$ is any graph with minimum degree $k$. We use the same construction, each $v\in V(H)$ independently chooses $k$ random $H$-neighbors to be placed in $H_{k-out}$. Thus if $H=K_n$ then $H_{k-out}=G_{k-out}$. Frieze and Johansson \cite{FrJo} proved that if $H$ has $n$ vertices and minimum degree at least $\brac{\frac12+\e}n$ then $H_{k-out}$ is Hamiltonian w.h.p. for $k\geq k_\e$. If $\e=0$ then two cliques of size $m$ intersecting in two vertices, shows that $k_0$ is not bounded as a function of $n$.
\begin{Problem}
Determine the growth rate of $k_0$. Suppose we assume also that $G$ has connectivity $\k\to \infty$. How fast should $\k$ grow so that $k_0=O(1)$.
\end{Problem}
Frieze, Karonski and Thoma \cite{FKT} considered the graphs induced by the unions of random spanning trees. They showed that 5 random trees are enough to guarantee Hamiltonicity w.h.p.
\begin{Problem}
Show that the union of 3 random spanning trees is enough to guarantee Hamiltonicity w.h.p.
\end{Problem}
Gao, Kami\'nski, MacRury and Pralat \cite{GKMP} consider the following model of a semi-random version of $G_{k-out}$. We start initially with the empty graph $G_0$. $G_{m}$ is obtained from $G_{m-1}$ as follows: $v_m\in [n]$ is chosen uniformly at random. Then a vertex $w_m$ is deterministically chosen and the edge $\set{v_m,w_m}$ is added to create $G_{m}$. In the paper they give a choice rule for $w_m$ that shows that w.h.p. $G_{2.61...n}$ is Hamiltonian. The paper \cite{BoF} has shown that $G_{3m}$ is Hamiltonian for a random choice of $w_m$. The results of \cite{GKMP} were improved by Frieze, Gao, MacRury, Pralat and Sorkin \cite{FGKMPS} to $\approx1.85n$. They also prove a lower bound of $\approx 1.26n$ for this problem. 
\begin{Problem}
Close the gap between the upper and lower bounds for this problem.
\end{Problem}
Frieze, Krivelevich and Michaeli \cite{builder} consider an on-line process for constructing various objects, Hamilton cycles in particular. They observe the graph process for up to $t$ steps and at each step they are required to decide irrevocably to accept or reject the next edege. There is a bound of $b$ on the number edges that can be accepted. A strategy for building an object in this model is called a $(t,b)$ strategy. They proved 
\begin{theorem}
  For every $\e>0$ there exists $C>1$ such that the following hold.
  \begin{enumerate}
    \item 
      If $t\ge (1+\e)n\log{n}/2$ and $b\ge Cn$ then there exists a $(t,b)$-strategy $B$ of Builder such that
      \[
        \lim_{n\to\infty} \Pr(B_t\text{ is Hamiltonian}) = 1.
      \]
    \item 
      If $t\ge Cn\log{n}/2$ and $b\ge(1+\e)n$ then there exists a $(t,b)$-strategy $B$ of Builder such that
      \[
        \lim_{n\to\infty} \Pr(B_t\text{ is Hamiltonian}) = 1.
      \]
\item Let $\t_2$ denote the hitting time for minimum degree two. There exists $C>1$ for which there exists a $(\tau_2,Cn)$-strategy $B$ of Builder with
  \[
    \lim_{n\to\infty} \Pr(B_{\tau_2}\text{ is Hamiltonian})=1.
  \]
  \end{enumerate}
\end{theorem}
Anastos \cite{AAA} showed that we can replace $C$ by $1+\e$ in 1.,2. of the above theorem.
\begin{Problem}
If $b=n+\om$, where $\om\to\infty$, how large should $t$ be in a $(t,b)$ strategy that constructs a Hamilton cycle w.h.p.
\end{Problem}
\subsection{The $n$-cube}
The graph $Q_{n}$ has been widely studied. Here $V(Q_n)=\set{0,1}^n$ and two vertices are adjacent if their Hamming distance is one. There are various models of random subgraphs of $Q_n$ and we mention two: in $Q_{n,p}^{(e)}$ we keep all the vertices of $Q_n$ and the edges of $Q_n$ independently with probability $p$. In $Q_{n,p}^{(v)}$ we choose a random subset of $V(Q_n)$, where each vertex is included independently with probability $p$. After this we take the subgraph induced by the chosen set of vertices. It is known for example that $Q_{n,p}^{(e)}$ becomes connected at around $p=1/2$. Also, Bollob\'as \cite{Bolmatch} determined the value of $p$ for there to be a perfect matching in $Q_{n,p}^{(e)}$ w.h.p., again at around $p=1/2$.  Condon, Espuny D\'iaz, Gir$\tilde{\text{a}}$o, K\"uhn and Osthus \cite{CDGKO2} proved that the threshold for the existence of a Hamilton cycle in  $Q_{n,p}^{(e)}$ is $p=1/2$. They verified property $\cA_k,k\geq 2$ and proved a perturbation result as in Section \ref{pert}.
\begin{Problem}\label{P30}
Determine the minimum value of $p$ for there to be a Hamilton cycle in $Q_{n,p}^{(v)}$, conditional on $Q_{n,p}^{(v)}$ containing as many odd as even vertices. Perhaps consider  resilience and color the edges.
\end{Problem}

\subsection{Random Lifts}
Amit and Linial \cite{AL} introduced the notion of a random lift of a fixed graph $H$. We let $A_v,v\in V(H)$ be a collection of sets of size $n$. Then for every $e=\set{x,y}\in E(H)$ we construct a random perfect matching $M_e$ between $A_x$ and $A_y$. The graph with vertex set $\bigcup_{v\in V(H)}A_v$ and edge set $\bigcup_{e\in E(H)}M_e$ is a {\em random lift} of $H$. 

Burgin, Chebolu, Cooper and Frieze \cite{BCCF} proved that if $s$ is sufficently large then a random lift of $K_s$ is Hamiltonian w.h.p. {\L}uczak, Witkowski and Witkowsi \cite{LWW} improved this and showed that if $H$ has minimum degree at least 5 and contains two edge disjoint Hamilton cycles then a random lift of $H$ is Hamiltonian w.h.p. This implies that a random lift of $K_5$ is Hamiltonian w.h.p. A random lift of $K_3$ consists of a set of vertex disjoint cycles.
\begin{Problem}\label{P31}
Is a random lift of $K_4$ Hamiltonian w.h.p.?
\end{Problem}
\subsection{Random Graphs from Random Walks}
Given a graph $G$, one can obtain a random set of edges by constructing a random walk. This was the view takne in Frieze, Krivelevich, Michaeli and Peled \cite{FKMP}. So, given $G$, we let $G_m$ denote the random subgraph of $G$ induced by the first $m$ steps of a simple random walk on $G$. The considered the case where $G=G_{n,p},p=\frac{C\log n}{n}$ and they showed that for every $\e>0$, there exists $C_\e$ such that $C\geq C_\e$ and $m\geq (1+\e)n\log n$ then w.h.p. $G_m$ is Hamiltonian. When $G=K_n$ they showed that w.h.p. $G_m$ is Hamiltonian for $m$ equal to one more than the number of steps needed to visit every vertex.
\begin{Problem}
Determine $C(\e)$ up to an $(1+o(1))$ factor.
\end{Problem}
\subsection{Random Geometric Graphs}
Let $X_1,X_2,\ldots,X_n$ be chosen independently and uniformly at random from the unit square $[0,1]^2$ and let $r$ be given. Let $\cX=\set{X_1,X_2,\ldots,X_n}$. The random geometric graph $G_{\cX,r}$ has vertex set $\cX$ and an edge $X_iX_j$ whenever $|X_i-X_j|\leq r$. See Penrose \cite{Pen} for more details or Chapter 11.2 of \cite{FK} for a gentle introduction. Di\'az, Mitsche and P\'erez-Gim\'enez \cite{DMP} showed that if $r\geq (1+\e)\bfrac{\log n}{\p n}^{1/2}$ then $G_{\cX,r}$ is Hamiltonian w.h.p. Balogh, Bollob\'as, Krivelevich, M\"uller and Walters \cite{BBKMW} proved that if we grow $r$ from zero then w.h.p. the ``hitting time'' for minimum degree at least two coincides with the hitting time for Hamiltonicity. M\"uller, P\'erez-Gim\'enez and Wormald \cite{MPW} proved that as $r$ grows, the hitting time for minimum degree $k$ coincides with the hitting time for property $\cA_k$, w.h.p. The papers \cite{BBKMW} and \cite{MPW} both deal with dimensions $d\geq 2$. The paper \cite{BBKMW} also deals with the nearest neighbor graph.
\begin{Problem}
Discuss resilience in the context of $G_{\cX,r}$.
\end{Problem}
Bal, Bennett, P\'erez-Gim\'enez and Pralat \cite{BBPP} considered the problem of the existence of a rainbow Hamilton cycle. They show that for $r$ at the threshold for Hamiltonicity, $O(n)$ random colors are sufficient to have a rainbow Hamilton cycle w.h.p. Frieze and P\'erez-Gim\'enez \cite{FP} reduced the number of colors to  $n+o(n)$.
\begin{Problem}
Consider the problem where there are exactly $n$ colors available.
\end{Problem}

Fountoulakis, Mitsche, M\"uller and Schepers \cite{FMMS} considered the KPKVB model. The points are chosen from a disk in the hyperbolic plane. The definition is somewhat complicated and can of course be found in \cite{FMMS}. There is a parameters $\a,\n$ and they prove that given $\a<1/2$ there are values $\n_0(\a),\n_1(\a)$ such that w.h.p. there is no Hamilton cycle if $\n<\n_0$ and there is a Hamilton cycle if $\n>\n_1$.
\begin{Problem}
Prove that $\n_0(\a)=\n_1(\a)$, as conjectured in \cite{FMMS}.
\end{Problem}
Ganeson \cite{Ganeson} considered random geometric graphs with $r$ just below the Hamiltonicity threshold and discussed the minimum number of edges needed to make the graph Hamiltonian.
\subsection{Random Intersection Graphs}
The random intersection graph $G_{n,m,p}$ is the intrersection graph of $S_1,S_2,\ldots,S_n$ where each $S_i$ is independently chosen as a subset of $[m]$ where an element is independently included with probability $p$. Hamiltonicity of $G_{n,m,p}$ and the related uniform model has been considered by Efthymiou and Spirakis \cite{ES}, Bloznelis and Radavi\v{c}ius \cite{BR} and by Rybarczyk \cite{Ryb1}, \cite{Ryb2}. In particular, \cite{Ryb1} proves
\begin{theorem}
Let $\a>1$ be constant and $m=n^\a$. Let $p_{\pm}=\sqrt{\frac{\log n+\log\log n\pm\om}{mn}}$ where $\om\to \infty$. Then w.h.p. $G_{n,m,p_-}$ is not Hamiltonian and w.h.p. $G_{n,m,p_+}$ is Hamiltonian.
\end{theorem}
Furthermore, in \cite{Ryb2} it is shown that w.h.p. the polynomial time algorithm HAM of \cite{BFF} is successful w.h.p. on $G_{n,m,p}$ whenever $m\gg \log n$ and $mp^2\leq 1$.
\begin{Problem}
Discuss resilience and other questions in relation to $G_{n,m,p}$.
\end{Problem}
\subsection{Preferential Attachment Graph}
The {\em Preferential Attachment Graph} (PAM) is a random graph sequence $G_0,G_,\ldots,G_n,\ldots,$ that bears some relation to networks found in the real world. Its main characteristic is having a heavy tail distribution for degrees. $G_{n+1}$ is obtained from $G_n$ by adding a new vertex $v_{n+1}$ and $m$ (a parameter) random edges. The distinguishing feature is that the $m$ neighbors of $v_{n+1}$ in $V(G_n)$ are chosen with probability proportional to their current degree. Frieze, Pralat, P\'erez-Gim\'enez and Reiniger \cite{FPPR} showed that if $m\geq 29500$ then $G_n$ is Hamiltonian w.h.p.
\begin{Problem}\label{P35}
Find the smallest $m$ such $G_m$ is Hamiltonian w.h.p.
\end{Problem}
\subsection{Nearest neighbor Graphs}
Given a graph $G=(V,E)$ we let $\s=(e_1,e_2,\ldots,e_m)$ be a random permutation of its edges. We think of the permutation being derived by giving each edge an independent random weight and then ordereing the edges in increasing order of weight. The random graph $G_{k-NN}$ is obtained as follows. For each $v\in V$ let $F_v=\set{f_1,f_2,\ldots,f_k}$ be the first $k$ edges in the sequence $\s$ that contain $v$. Let $F=\bigcup_vF_v$ and then $G_{k-NN}=(V,F)$. When $G=K_n$, Cooper and Frieze \cite{CF-NN} determined the connectivity and when $G$ is a random geometric graph, Balister, Bollob\'as, Sarkar and Walters \cite{BBSW} determined the connectivity approximately. (Exact determination in this case is a nice open problem.)
\begin{Problem}
Determine, for various $G$, the minimum $k$ for which $G_{k-NN}$ is Hamiltonian. When $G=K_n$ this should be 3 or 4. When $G$ is a random geometric graph, this should be $c\log n$. (This is likely to be very difficult, seeing as the connectivity threshold is still open.)
\end{Problem}
We noe that \cite{builder} shows that the minimum is $O(1)$ for the case $G=K_n$.
\subsection{Achlioptas Process}
In this model, sets of $K$ random edges are presented sequentially and one is allowed to choose one in order to fulfill some purpose. Call each choice a round. Krivelevich, Lubetzky and Sudakov \cite{KLSA} considered the problem of optimizing the selection so that one can obtain a Hamilton cycle as quickly as possible. They show (i) if $K\gg\log n$ then $n+o(n)$ rounds are sufficient w.h.p. and (ii) if $K=\g\log n$ then w.h.p. the number of rounds $\t_H$ satisfies
\beq{Ach}{
1+\frac{1}{2\g}+o(1)\leq \frac{\t_H}{n}\leq 3+\frac{1}{\g}+o(1).
}
Anastos \cite{AAA} showed that $\t_H/n\sim 1+1/2\g$ w.h.p.
\subsection{Maximum degree process}
In the maximum degree $d$ process, edges are added to the empty graph on vertex set $[n]$, avoiding adding edges that make the maximum degree more than $d$. For $d=2$, Telcs, Wormald and Zhou \cite{TWZ} showed that the probability the process terminates with a hamilton cycle is asymptotically equal to $c_1n^{-1/2}$ for an explicitly defined $c_1$.
\subsection{Pancyclicity}
A graph with $n$ vertices is {\em pancyclic} if it contains cycles of lengths $3\leq k\leq n$. Cooper and Frieze \cite{CF4} showed that the limiting probability for $G_m$ to be pancyclic is the same as the limiting  for minimum degree at least two.  This was refined by Cooper \cite{C1}. He showed that w.h.p. there is a Hamilton cycle $H$ such that cycles of every length can be constructed out of the edges of $H$ and at most two other edges per cycle. Then in \cite{C1a} Cooper showed that one edge per cycle is suffcient. Lee and Samotij \cite{LS1} determined the resilience of pancyclicity. They show that if $p\geq n^{-1/2}$ then w.h.p. every Hamiltonian subgraph $G'\subseteq G_{n,p}$ with more than $(1/2+o(1))n^2p/2$ edges is pancyclic.
\begin{Problem}
Determine the threshold for $G_{n,p}$ to contain the $r$th power, ($r\geq 2$), of a cycle of length $k$ for all $2\leq k\leq n$. 
\end{Problem}
Krivelevich, Lee and Sudakov \cite{KLS2} proved that $G=G_{n,p},p\gg n^{-1/2}$ remains pancyclic w.h.p. if a subgraph $H$ of maximum degree $(\frac12-\e)np$ is deleted, i.e. pancyclicity is locally resilient. The same is true for random regular graphs when $r\gg n^{1/2}$.

Alon and Krivelevich \cite{AKpan} proved that if $p\geq (1+o(1))\log n/n$ then w.h.p. $G_{n,p}$ contains a pancyclic subgraph with $n+(1+o(1))\log_2n$ edges. In \cite{AKpan1} they prove that w.h.p. $\m(G_{n,p})=\wh\m(G_{n,p}$ if $np\geq 20$. Here $\m(G)$ is the minimum number of edges additional edges needed to make $G$ Hamiltonian and $\wh\m$ is the corresponding number to make $G$ pancyclic. They also show that if $np=d$, $\m(G_{n,p})\sim f(d)n$ where $f(d)=de^{-d}/2+\cdots$.
\subsection{Hamilton Game}
There is a striking and mysterious relationship between the existence of Hamilton cycles and the Hamilton Maker-Breaker game. In this game played on some Hamiltonian graph $G$, two players Maker and Breaker take turns in selecting (sets of) edges. Maker tries to obtain the edges of A Hamiltonian subgraph and Breaker tries to prevent this.  There is a bias $b$ for breaker, if Breaker is allowed to choose $b$ edges for every choice by Maker. Ben-Shimon, Ferber, Hefetz and Krivelevich \cite{BFHK12} prove a hitting time result for the $b=1$ Hamilton cycle game on the graph process. Assuming that Breaker starts first, Maker will have a winning strategy in $G_m$ iff $m\geq m_4$, the hitting time for minimum degree 4. This is best possible. Biased Hamiltonicity games on $G=G_{n,p}$ were considered in Ferber, Glebov, Krivelevich and Naor \cite{FGKN15}  where it was shown that for $p\gg\frac{\log n}{n}$, the threshold bias $b_{HAM}$ satisfies $b_{HAM}\approx \frac{np}{\log n}$ w.h.p. 

Hefetz, Krivelevich and Tan \cite{HKT} considered a variant on this game. In the $(1:q)$ Waiter-Client version, in each round, Waiter offers Client $q+1$ previously unoffered edges and Clint chooses one. Waiter wins if he can force Client to choose a Hamiltonian graph.  Let $\cW_q$ denote the property that there is a winning strategy for waiter. This is a monotone property and they show that $\frac{\log n}{n}$ is a sharp threshold for this property when the game is played on $G_{n,p}$. In the Client-Waiter game, Client wins if he can claim a Hamilton cycle. In this game $\frac{(q+1)\log n}{n}$ is a sharp threshold.
\section{Random Digraphs $D_{n,m}$ and $D_{n,p}$}
The random graphs $D_{n,m}$ and $D_{n,p}$ are as one might expect, directed versions of $G_{n,m},G_{n,p}$ repectively. For $D_{n,m}$ we choose $m$ random edges from the complete digraph $\vec{K}_n$ and for $D_{n,p}$ we include each of the $n(n-1)$ edges of $\vec{K}_n$ independently with probability $p$. 
\subsection{Existence}
The existence question was first addressed by Angluin and valiant \cite{AV}. They showed that if $p\geq \frac{K\log n}{n}$ for sufficiently large $K$ then $D_{n,p}$ is Hamiltonian w.h.p. Using an elegant interpolation between $D_{n,p}$ and $G_{n,p}$, McDiarmid \cite{McD0} proved that for any $0\leq p\leq 1$,
\[
\Pr(D_{n,p}\text{ is Hamiltonian})\geq \Pr(G_{n,p}\text{ is Hamiltonian}).
\]
This shows that if $p=\frac{\log n+\log\log n+\om}{n}$ then $D_{n,p}$ is Hamiltonian w.h.p. Frieze \cite{Fdir} proved a hitting time result. Consider the directed graph process $D_0,D_1,\ldots,D_{n(n-1)}$ where $D_{m+1}$ is obtained from $D_m$ by adding a random directed edge. Let $m_{\cH}$ be the minimum $m$ such that $D_m$ is Hamiltonian and let $m_k$ be the minimum $m$ such that $\delta^{\pm}(D_m)\geq k$ where $\delta^+$ and $\delta^-$ denote minimum out- and in-degree respectively. Then \cite{Fdir} shows that w.h.p. $m_{\cH}=m_1$ w.h.p. This removes a $\log\log n$ term from McDiarmid's result.

Cooper and Frieze \cite{CFDham} proved that if $\cdnm{1}$ denotes the set of digraphs with vertex set $[n]$ and $m$ edges such that $\min\set{\delta^-,\delta^+}\geq 1$. Then let $\dnm{1}$ be sampled uniformly from $\cdnm{1}$.
\begin{theorem}\label{th1}
Let $m=\tfrac n2(\log n+2\log\log n+c_n)$ then
\beq{H2}{
\lim_{n\to\infty}\Pr(\dnm{1}\text{ is Hamiltonian})=\begin{cases}0&c_n\to-\infty.\\e^{-e^{-c}/4}&c_n\to c.\\1&c_n\to\infty.\end{cases}
}
\end{theorem}
The R.H.S. of \eqref{H2} is the limiting probability that $\dnm{1}$ contains two vertices of in-degree one (resp. out-degree one) that share a common in-neighbour (resp. common out-neighbor).
Ferber, Sales and Shurman \cite{FSS} considered the problem of covering the edges of $D_{n,p}$ in Hamilton cycles. They show that if $1/2\geq p\geq \frac{\log^{20}n}{n}$ then w.h.p. this can be done using the maximum of the in-degree and out-degree of $D_{n,p}$.
\begin{Problem}
Tighten the results of \cite{FSS}. Is there an optimal covering at the hittng time for Hamiltonicity?
\end{Problem}
\subsection{Packing, Covering and Counting}
The paper \cite{Fdir} shows that w.h.p. at time $m_k$, $D_m$ contains $k$ edge disjoint Hamilton cycles. Ferber, Kronenberg and Long \cite{FKL} proved that if $np/\log^4n\to\infty$ then w.h.p. $D_{n,p}$ contains $(1-o(1))np$ edge disjoint Hamilton cycles. This was improved by Ferber and Long \cite{FeLo}, see below.
\begin{Problem}
Let $\delta=\min\set{\delta^+,\delta^-}$. Is it true that throughtout the directed random graph process $D_m,m\geq 0$ that w.h.p. $D_m$ contains $\delta$ edge disjoint Hamilton cycles?
\end{Problem}
The paper \cite{FKL} also considered covering the edges of $D_{n,p}$ by Hamilton cycles. They show that if $p\gg\frac{\log^2n}{n}$ then the edges of $D_{n,p}$ can be covered by $(1+o(1))np$ Hamilton cycles. The paper by Ferber, Sales and Sherman \cite{FSS} showed that $D_{n,p}$ can be covered with $\delta$ Hamilton cycles w.h.p., provided that  $np\geq \log^{20}n$.
\begin{Problem}
Is it true that if $p\gg\frac{\log n}{n}$ then the edges of $D_{n,p}$ can be covered by $(1+o(1))np$ Hamilton cycles.
\end{Problem}
Finally, consider the number of Hamilton cycles in $D_{n,p}$. The paper \cite{FKL} shows that if $p\gg\frac{\log^2n}{n}$ then w.h.p. $D_{n,p}$ contains $(1+o(1))^nn!p^n$ Hamilton cycles. Ferber, Kwan and Sudakov \cite{FKS} improved this to show that w.h.p. at the hitting time for the existence of a directed Hamilton cycle, there are w.h.p. $(1+o(1))^nn!p^n$ disitinct Hamilton cycles.
\begin{Problem}
At the hitting time for the existence of a Hamilton cycle, $D_m$ w.h.p. contains $\a_nn!p^n$ Hamilton cycles. Determine $\a_n$ as accurately as possible.
\end{Problem}
Ferber and Long \cite{FeLo} considered Hamilton cycles with arbitrary orientations of the edges. They showed that if $C_1,C_2,\ldots,C_t,t\leq (1-\e)np$ are arbitrarily oriented Hamilton cycles  and if  $np/\log^3n\to\infty$ then w.h.p.$D_{n,p}$ contains edges disjoint copies of these cycles. They also show that w.h.p. $D_{n,p}$ contains $(1+o(1))^nn!p^n$ copies of any arbitrariliy oriented cycle. They conjectured the truth of the following:
if $np-\log n\to \infty$ and $C$ is some arbitrarily oriented Hamilton cycle, then $D_{n,p}$ contains a copy of $C$ w.h.p. Frieze, Pralat and P\'erez-Gim\'enez \cite{FPP} studied the existence of Hamilton cycles where the orientation of edges follows a repeating pattern. They proved hitting time versions. In particular, for the pattern where the orientations alternate, they show that approximately $\frac{1}{2}n\log n$ random edges are needed. 

Montgomery \cite{M3} strengthened \cite{FPP} considerably and proved the conjecture of \cite{FeLo}.  In addition he showed that w.h.p. the last pattern to appear is that of a properly oriented cycle i.e. one in which every vertex has in- and out-degree one. Araujo, Balogh, Krueger, Piga and Treglown \cite{ABKPT} considered the following model: start with a digraph $D_0$ that  is an orientation of an $n$ vertex simple graph $G$ such that the minimum in- and out-degree is at least $\a n$ and add $C(\a)n$ random directed edges to create a digraph $D$. They show that w.h.p. $D$ has a Hamilton cycle of every possible orientation.

Gishboliner, Krivelevich and P. Michaeli \cite{GKM1} study a form of discrepancy related ot Hamilton cycles. In particular they proved that if $p$ is above that Hamiltonicity threshold then  w.h.p. in any orientation of $G_{n,p}$ there is a Hamilton cycle with almost all edges oriented in the same direction. 

One can also consider cores in the context of digraphs. The $k$-core of a digraph $D$ will be the largest subgraph with in- and out-degree at least $k$.
\begin{Problem}
For which values of $k$ are the cores of $D_{m}$ born Hamiltonian.
\end{Problem}
\subsection{Resilience}
Hefetz, Steger and Sudakov \cite{HSS} began the study of the resilience of Hamiltonicity for random digraphs. They showed that if $p\gg\frac{\log n}{n^{1/2}}$ then w.h.p. the Hamiltonicity of $D_{n,p}$ is resilient to the deletion of up to $(\frac12-o(1))np$ edges incident with each vertex. The value of $p$ was reduced to $p\gg\frac{\log^8n}{n}$ by Ferber, Nenadov, Noever, Peter and \v{S}kori\'c \cite{FNNPS}. Finally, Montgomery \cite{M2} proved that in the random digraph process, at the hitting time for Hamiltonicity, the property is resilient w.h.p.
\section{Other models of Random Digraphs}
\subsection{$k$-in,$k$-out}
The random graph $D_{k-in,\ell-out}$ is generated as follows. Each $v\in [n]$ independently chooses $k$ in-neighbors and $\ell$ out-neighbors. It is a directed version of the model $G_{k-out}$ considered in Section \ref{kout}. Cooper and Frieze \cite{CF3inout} showed that $D_{3-in,3-out}$ is Hamiltonian w.h.p. And then in \cite{CF2inout} they  showed that $D_{2-in,2-out}$ is Hamiltonian w.h.p. This is best possible, since w.h.p. $D_{1-in,1-out}$ is not Hamiltonian.
\begin{Problem}
The proofs in \cite{CF3inout}, \cite{CF2inout} can be seen as the analysis of an $n^{O(\log n)}$ time algorithm. Is there a polynomial time algorithmic proof?
\end{Problem}
The related directed nearets-neighbor digraph is relatively unexplored, although \cite{BBSW} does consider a directed version).
\subsection{Random Regular Digraphs}
Cooper, Frieze and Molloy \cite{CFM} proved that w.h.p. the random regular digraph $D_{n,r}$ is hamiltonian for every fixed $r\geq 3$. In $D_{n,r}$ each vertex $v\in [n]$ has in-degree and out-degree $r$. 
\begin{Problem}
Discuss the Hamiltonicity of $D_{n,r}$ for $r\to\infty$ with $n$.
\end{Problem}
\begin{Problem}
Discuss the query complexity, (in the context of \cite{FKSV}), of finding Hamilton cycles in the random digraph $D_{n,p}$. Is it $n+o(n)$?
\end{Problem}
Another way to generate random regular digraphs, is to take the union of $r$ random permuation digraphs. \cite{Fperm} shows that the union of 3 directed permutation digraphs is Hamiltonian. Cooper \cite{Cperm} showed that 2 is not enough.
\subsection{$D_p$}
In the same way that we defined $G_p$ a a random subgraph of an arbitrary graph $G$, we can define $D_p$ as a random subgraph of an arbitrary digraph $D$. In particular, similarly to Problem \ref{Johan}, we can pose
\begin{Problem}
Let $D$ be a digraph with vertex set $[n]$ and minimum out- and in-degree at least $(\frac12+\e)n$. Now consider the random digraph process restricted to the edges of $D$. Is the hitting time for Hamiltonicity equal to the hitting time for out- and in-degree at least one, w.h.p.?
\end{Problem}
\subsection{Hamilton Game}
Frieze and Pegden \cite{FPgame} considered the Maker-Breaker game on the complete digraph $\vec{K}_n$. They showed that if $b\geq \frac{(1+\e)\log n}{n}$ then Breaker wins and that if $b\leq \frac{\e\log n}{n}$ then Maker wins, $\e$ sufficiently small.
\begin{Problem}
Show that Maker wins if $b\leq \frac{(1-\e)\log n}{n}$
\end{Problem}
\subsection{Random Lifts}
Given a digraph $D=(V,E)$ we can construct a random lift as follows: We let $A_v,v\in V$ be a collection of sets of size $n$. Then for every oriented edge $e=(x,y)\in E(H)$ we construct a random perfect matching $M_e$ between $A_x$ and $A_y$. The edges of this matching are oriented from $A_x$ to $A_y$. Chebolu and Frieze \cite{Chef} proved that if $H=\vec{K}_h$ for a sufficiently large $h$, then a random lift of $H$ is Hamiltonian w.h.p.
\begin{Problem}
Show that a random lift of $\vec{K}_3$ is Hamiltonian w.h.p.
\end{Problem} 
\subsection{Random Tournaments}
K\"uhn and Osthus \cite{KO1} showed a random tournament contains $\delta$ edge disjoint Hamilton cycles, where $\delta=\min\set{\delta^+,\delta^-}$ and $\delta^+$ denotes the minimum out- and $\delta^-$ denotes the minimum in-degree.
\subsection{Perturbations of Dense Digraphs}
Krivelevich, Kwan and Sudakov \cite{KKS} show that if $D$ is a digraph with vertex set $[n]$ and minimum in- and out-degree at least $\a n$ and $R$ is a set of $c=c(\a)$ random directed edges, then w.h.p. $D+R$ is Hamiltonian, indeed, pancyclic. They also consider random perturbations of a tournament. Suppose that $T$ is a tournament on vertex set $[n]$ in which each in- and out-degree is at least $d$. Now independently choose $m\gg \frac{n}{d+1}$ random edges of $T$ and then orient them uniformly at random. Then w.h.p. the resulting perturbed tournament has at least $q$-edge disjoint Hamilton cycles, for any fixed $q$.

Araujo, Balogh, Krueger, Piga and Treglown \cite{ABKPT1} proved that in the same model as \cite{KKS}, the digraph $D+R$ simultaneously contains Hamilton cycles of every orientation. In addition they prove the same for cycles of length $k\in [2,(1-\eta)n]$.
\section{Edge-colored digraphs}
Anastos and Briggs \cite{AnBr} extended the result of \cite{BFKLS} to $D_{n,p}$. I.e. they give an on-line algorithm for constructing $k$ edges disjoint Hamilton cycles in the random digraph process at the point where the minimum in- and out-degrees are both first at least $k$.

Katsamaktis, Lezter and Sgueglia \cite{KaLS1} proved that if $D$ is an edge-colored random permutation of a dense digraph then $D$ contains a directed rainbow Hamilton cycle w.h.p.
\begin{Problem}
When do random digraphs contain rainbow copies of Hamilton cycles with arbitrary given orientations?
\end{Problem}
\section{The random hypergraphs $H_{n,m:k}$ and $H_{n,p:k}$}
In the main when consideringnhypergraphs, we will consider random $k$-uniform hypergraphs where each edge has size $k\geq 3$. The random hypergraph $H_{n,m:k}$ has vertex set $[n]$ and $m$ randomly chosen edges from $\binom{[n]}{k}$. Similiarly, the random hypergraph $H_{n,p:k}$ has vertex set $[n]$ and each element of $\binom{[n]}{k}$ is included independently as an edge with probability $p$. When $m=\binom{n}{k}p$, the two models behave similarly.

Suppose that $1\leq \ell< k$. An {\em $\ell$-overlapping Hamilton cycle} $C$in a $k$-uniform hypergraph $H=(V,\cE)$ on $n$ vertices is acollection of $m_\ell=n/(k-\ell)$ edges of $H$ such that for some cyclic order of $[n]$ every edge consists of $k$ consecutive vertices and for every pair of consecutive edges $E_{i-1},E_i$ in $C$ (in the natural ordering of the edges) we have $|E_{i-1}\cap E_i|=\ell$. Thus, in every
$\ell$-overlapping Hamilton cycle the sets $C_i=E_i\setminus E_{i-1},\,i=1,2,\ldots,m_\ell$, are a partition of $V$ into sets of size $k-\ell$. Hence, $m_{\ell}=n/(k-\ell)$. We thus always assume, when discussing $\ell$-overlapping Hamilton cycles, that this necessary condition, $k-\ell$ divides $n$,  is fulfilled. In the
literature, when $\ell=k-1$ we have a {\em tight} Hamilton cycle and when $\ell=1$ we have a {\em loose} Hamilton cycle.
\subsection{Existence}
Frieze \cite{Fhyp1} showed that if $K$ is sufficiently large and $4)|n$ then w.h.p. $H_{n,Kn\log n:3}$ contains a loose Hamilton cycle. Dudek and Frieze \cite{DuFr1} generalised the argument of \cite{Fhyp1} and showed that if $K$ is sufficiently large and $2(k-1)|n$ then w.h.p. $H_{n,Kn\log n:k}$ contains a loose Hamilton cycle. The divisibility conditions in these papers are not optimal and Dudek, Frieze, Loh and Speiss \cite{DFLS} relaxed these conditions to $(k-1)|n$.

Dudek and Frieze \cite{DuFr2} found the existence thresholds for all integers $k> \ell\geq 2$ up to a constant factor (except for $\ell=2$).  Narayanan and Schact \cite{NaSc} tightened these results further and proved the following: let $k>\ell>1$ and $s=k-\ell$ and $t=k-\ell\mod s$ and  $\l(k,\ell)=t!(s-t)!$. Let 
\[
p^*_{k,\ell}(n)=\frac{\l(k,\ell)e^s}{n^s}.
\]
Then, for any fixed $\e>0$,
\[
\Pr(H_{n,p;k}\text{ contains an $\ell$-overlapping Hamilton cycle})\to\begin{cases}1&p\geq(1+\e)p^*_{k,\ell}(n). \\0&p\leq(1-\e)p^*_{k,\ell}(n).\end{cases}
\]
The lower bound is proved by the first moment method and the upper bound is via a clever refinement of the second moment method. Frieze and P\'erez-Gim\'enez \cite{FP1} showed that the threshold for the existence of a loose Hamilton cycle in $H_{n,m:k}$ is asymptotically equal to $\frac{1}{r}n\log n$.
\begin{Problem}
Tighten the statements on the existence of Hamilton cycles. For the case $\ell=1$, prove a hitting time result.
\end{Problem}

\begin{Problem}
Determine the resilience of Hamiltonicity in random hypergraphs.
\end{Problem}
The minimum $d$-degree $\delta_d(H)$ of a hypergraph $H$ is equal to $\min_{S\subseteq V(H), |S|=d}|\set{e\in E(H):e\supseteq S}|$. The minimum $d$-degree threshold for loose Hamilton cycles $\m_d(k)$ is defined as the least $\m\in [0, 1]$ such that for every $\g > 0$ and large enough $n$ divisible by $k - 1$, every $n$-vertex $k$-unifom hypergraph $G$ with $\delta_d(G)\geq (\m + \g)\binom{n-d}{k-d}$ contains a loose Hamilton cycle. Alvarado, Kohayakawa, Lang, Mota and Stagni \cite{AKLMS} proved that or every $1 \leq d < k$ and $\g > 0$ there is a $C > 0$ such that the following holds. If $p \geq \max\{n^{-(k-1)/2+\g} , Cn^{-(k-d)} \log n\}$ and $n$ is divisible by $k - 1$, then w.h.p. $G=H_{n, p;k}$ has the property that every spanning subgraph $G' \subseteq  G$ with $\delta_d(G') \geq (\m_d(k) +\g)p\binom{n-d}{k-d}$ contains a loose Hamilton cycle.

Let $H_{n,m:k}^{(\ell)}$ denote a random $k$-uniform hypergraph with vertex set $[n]$,  $m$ edges and minimum degree at least $\ell$.
\begin{Problem}\label{P49}
Show that $H_{n,cn:k}^{(3)}$ is Hamiltonian w.h.p. for large enough $c$.
\end{Problem}
Parczyk and Person \cite{PaPe} proved an extension of Riordan's spanning subgraph result \cite{R} to hypergraphs. Among other things this gives a slightly weaker version of the known results on the thresholds for Hamilton cycles (except for loose) in random uniform hypergraphs and also applies to powers of tight Hamilton cycles. Chang, Han and Sun \cite{CHS} improved a result in \cite{PaPe} by showing that if $p\geq Cn^{-1/\binom{k+r-2}{k-1}}$ for some contant $C>0$, then w.h.p. $H_{n,p;k}$ contains the $r$th power of a tight Hamilton cycle.

Petrova and Truji\'c \cite{PeTr} proved the following for $d\in\set{1,2}$. For every $\g>0$ there is a $C>0$ such that the following holds: Suppose that $p\geq C\log n\cdot\max\set{n^{-3/2},n^{-3+d}}$ and $n$ is even. Then w.h.p. $H=H_{n,p:3}$ has the property that every spanning subgraph $G$ of $H$ with minimum $d$-degree $\delta_d(G)\geq (\delta_d(H)+\g)p\binom{n-d}{3-d}$ contains a loose Hamilton cycle. Here $\delta_d(H)$ is the minimum over sets $S$ of size $d$ of the number of edges containing $S$.
\subsection{$H_p$}
Here we consider random subgraphs of dense hypergraphs. For $1\leq d\leq k-1$ we let $\delta_d$ be the largest integer $m$ such that every set of $d$ vertices is contained in $m$ edges of $H$. It is shown in \cite{KMP} and \cite{JLS} that if $\delta_{k-1}(H)\geq (1/2+\e)n$ and if $p\gg 1/n$ then $H_p$ contains a tight Hamilton cycle w.h.p. The paper \cite{JLS} contains several more interesting results along these lines, including results on loose Hamilton cycles.
\subsection{Algorithms}
Allen, B\"ottcher, Kohayakawa and Person \cite{ABKP} gave a randomised polynomial time algorithm for finding a tight Hamilton cycle in $H_{n,p:k}$ provided $p\geq n^{-1+\e}$ for any fixed $\e>0$. Allen, Koch, Parczyk and Person \cite{AKPP} gave a deterministic polynomial time algorithm for finding a tight Hamilton cycle provided $p\geq \frac{C\log^3n}{n}$ for sufficiently large $C$.
\begin{Problem}
Construct a polynomial time algorithm for finding Hamilton cycles in random hypergraphs, for all relevant $\ell$ and $p$. 
\end{Problem}
Molloy, Pralat and Sorkin \cite{MPS} consider the following semi-random process. The player is trying to build an $s$-uniform hypergraph that contains a perfect matching or a Hamilton cycle. Each edge consists of $r$ randomly chosen vertices plus $s-r$ vertices chosen by the player. Among their results, they prove that if $1\leq r\leq2,s\geq3$ then only $O(n)$ edges are needed w.h.p. 
\begin{Problem}
Suppose that $s\geq 4$. Can we succeed w.h.p. in $O(n)$ steps if $r=s-1$?
\end{Problem}
\subsection{Random regular Hypergraphs}
Altman, Greenhill, Isaev and Ramadurai \cite{AGIR} determined the threshold degree for a random $r$-regular $k$-uniform hypergraph $H_{n,r:k}$ to have a  loose Hamilton cycle. In this paper, $r=O(1)$ and they prove
\[
\lim_{n\to\infty}\Pr(H_{n,r:k}\text{ contains a loose Hamilton cycle})=\begin{cases}1&r>\r(k).\\0&r\leq \r(k).\end{cases}
\]
Here $\r=\r(k)$ is the unique real in $(2,\infty)$ such that
\[
(\r-1)(k-1)\bfrac{\r k-r-k}{\r k-\r}^{k-1)(\r k-\r-k)/k}=1.
\]
Dudek, Frieze, Ruci\'nski and \v{S}ileikis \cite{DFRS} show that if $n\log n \ll m\ll n^k$ and $r\approx km/n$ then there is an embedding of $G_{n,m:k}$ into $H_{n,r:k}$ showing the existence of Hamilton cycles in $G_{n,r:k}$ w.h.p. whenever there is one w.h.p. for the corresponding $G_{n,m:k}$. Espuny D\'iaz, Joos, K\"uhn and Osthus \cite{DJKO} proved that if $2\leq \ell<k$ and $r\ll n^{\ell-1}$ then w.h.p. $H_{n,r:k}$ does not contain an $\ell$-overlapping Hamilton cycle.
\subsection{Rainbow Hamilton Cycles}
Let $H_{n,p,:k}^{(r)}$ be $H_{n,p:k}$ with its edges randomly colored from $[r],r=cn\geq 1/(k-\ell)]$. Ferber and Krivelevich \cite{FeKr} proved the following: Let $k>\ell \ge 1$ be integers. Suppose that $n$ is a multiple of $k-\ell$. Let $p\in [0,1]$ be such that w.h.p.  $H_{n,p:k}$ contains an $\ell$-overlapping Hamilton cycle. Then, for every $\varepsilon = \varepsilon(n)\ge 0$, letting $r= (1+\varepsilon) m_\ell$ and $q = r p / (\varepsilon m_\ell + 1)$ we have that w.h.p. $H_{n,q,:k}^{(cn)}$  contains a rainbow $\ell$-overlapping Hamilton cycle. 

This was improved by Dudek, English and Frieze \cite{DEF} to the following: \\
(i) Let $k > \ell \ge 2$ and $\e>0$ be fixed: (i) for all integers $k> \ell\geq 2$, if
\[
p\leq
\begin{cases}
 (1-\e)e^{k-\ell+1}/n^{k-\ell} & \text{ if } c = 1/(k-\ell)\\
 (1-\e)\bfrac{c-1/(k-\ell)}{c}^{(k-\ell)c-1}e^{k-\ell+1}/n^{k-\ell} & \text{ if } c > 1/(k-\ell),
\end{cases}
\] 
then w.h.p. $H_{n,p,k}^{(cn)}$ is not rainbow $\ell$-Hamiltonian.\\
(ii) For all integers $k>\ell \ge 3$, there exists a constant $K=K(k)$ such that if $p\geq K/n^{k-\ell}$
and $n$ is a multiple of $k-\ell$ then  $H_{n,p:k}^{(cn)}$ is rainbow $\ell$-Hamiltonian w.h.p.\\
(iii) If $k>\ell=2$ and $n^{k-1}p\to\infty$ and $n$ is a multiple of $k-2$, then $H_{n,p:k}^{(cn)}$ is rainbow $2$-Hamiltonian w.h.p.  \\
(iv) For all $k\geq 4$, if
\[
p\geq 
\begin{cases}
(1+\e)e^2/n & \text{ if } c=1\\
(1+\e) \bfrac{c-1}c^{c-1} e^2 / n & \text{ if } c>1,
\end{cases}
\] 
then w.h.p. $H_{n,p,k}^{(cn)}$ is rainbow $(k-1)$-Hamiltonian, i.e. it contains a rainbow \emph{tight} Hamilton cycle.\\
(v) Fix $k\ge 3$ and suppose that $(k-1)|n$. Let $r=n/(k-1)$ and $n^{k-1}p/\log n\to\infty$. Then, w.h.p. $H_{n,p:k}^{(cn)}$ contains a rainbow loose Hamilton cycle.
\begin{Problem}
How large should $p$ be, so that w.h.p. $H_{n,p:k}^{m_\ell}$ contains an $\ell$-overlapping Hamilton cycle.
\end{Problem}
\subsection{Perturbations of dense hypergraphs}
In this section we consider adding random edges to suitably dense hypergraphs.  McDowell and Mycroft \cite{McM} proved that for integers $2\leq \ell<k$ and a small constant $c$, the union of a $k$-uniform hypergraph with linear minimum codegree and $H_{n,p:k},p\geq n^{-(k-\ell-c)}$ contains an n $\ell$-overlapping Hamilton cycle w.h.p. Bedenknecht, Han, Kohayakawa and Mota \cite{BHKM} proved the following: For $k\geq 2$ and $r\geq 1$ such that $k+r\geq 4$, and for any $\a>0$, there exists $\e>0$ such that the union of an $n$-vertex $k$-uniform hypergraph with minimum codegree $(1-(k+r-2k-1)-1+\a)n$  and $G_{n,p:k}$ with $p\geq n-(k+r-2k-1)-1-\e$ on the same vertex set contains the $r$th power of a tight Hamilton cycle w.h.p.
Chang, Han and Thoma \cite{CHT} extended the result of \cite{BHKM} and proved that for $k\geq 3,r\geq 2$ and $\a>0$ there exists $\e$ such that the following holds: suppose that $H$ is a $k$-uniform hypergraph on $n$ vertices such that every set of $(k-1)$ vertices is contained in at least $\a n$ edges and $p\geq n^{-\binom{k+r-2}{k-1}^{-1}-\e}$ then w.h.p. $H+H_{n,p:k}$ contains the $r$th power of a tight Hamilton cycle.

Krivelevich, Kwan and Sudakov \cite{KKS} proved that if the $k$-uniform hypergraph $H$ is such that every set of $(k-1)$ vertices is contained in at least $\a n$ edges then there exists $c_k=c_k(\a)$ such that if $R$ consists of $c_kn$ random edges, then w.h.p. $H+R$ contains a loose Hamilton cycle. 
\subsection{Other types of Hamilton cycle}
A {\em weak Berge Hamilton cycle} is a sequence $v_1,e_1,v_2,\ldots,v_n,e_n$ of vertices $v_1,v_2,\ldots,v_n$ where $v_1,v_2,\ldots,v_n$ is a permutation of $[n]$ and $e_1,e_2,\ldots,e_n$ are edges such that $e_i$ contains $\set{v_i,v_{i+1}}$. We drop ``weak'' if the edges are distinct. Poole \cite{Poole} proved that if $p=(k-1)!\frac{\log n+c_n}{n^{r-1}}$ then
\[
\lim_{n\to\infty}\Pr(H_{n,p:k}\text{ contains a weak Berge Hamilton cycle})=\begin{cases}0&c_n\to-\infty.\\ e^{-e^{-c}}&c_n\to c.\\1&c_n\to\infty.\end{cases}
\]
Bal and Devlin \cite{BalD} were off by a factor of $k$ for the threshold for Berge Hamilton cycles and the exact threshold was settled by Bal, Berkowitz, Devlin and Schacht \cite{BBDS} who proved a hitting version for both Berge Hamilton cycles and weak Berge Hamilton cycles.

Bal and Devlin and Bal, Berkowitz, Devlin and Schacht also considered the random hypergraph  $H_{r-out}$. Here each vertex $v$ randomly chooses $r$ edges containing $v$. They showed that if $k\geq 4$ then $H_{r-out}$ contains a Berge Hamilton cycle w.h.p. if and only if $r\geq 2$. They also show that if $k\geq 3$ then  $H_{r-out}$ contains a weak Berge Hamilton cycle w.h.p. if and only if $r\geq 2$.
\begin{Problem}
Prove packing and covering results for (weak) Berge Hamilton cycles.
\end{Problem}
\begin{Problem}
Prove algorithmic versions of the results in \cite{BBDS}.
\end{Problem}
\begin{Problem}
Answer resilience questions relating to the results of \cite{BBDS}.
\end{Problem}

Clemens, Ehrenm\"uller and Person \cite{CEP} proved a Dirac type of result. Suppose that $k\geq 3,\,\g>0$ and $p\geq \frac{\log^{17r}n}{n^{r-1}}$. Let $H$ be a spanning subgraph of $H_{n,p:k}$ with minimum vertex degree at least $\brac{\frac{1}{2^{k-1}}+\g}\binom{n}{k-1}p$. Then w.h.p. $H$ contains a Berge Hamilton cycle. The minimum degree is tight in the sense that one cannot replace the $+\g$ by $-\g$ for some small $\g$. Im and Kim \cite{ImKim} proved that if an $r$-uniform hypergraph $H$ has minimum degree at least $\brac{\frac{1}{2^{r-1}}+\e}\binom{n-1}{r-1}$ then the hitting time for minimum degree 2 coincides w.h.p. with the hitting time for a Berge Hamilton cycle.  
\begin{Problem}
Optimize the $\log^{O(1)}n$ factor in \cite{CEP}.
\end{Problem}

Dudek and Helenius \cite{DuH} considered {\em offset Hamilton cycles}. An $\ell$-offset hamilton cycle in a $k$-uniform hypergraph is a sequence of edges $E_1,E_2,\ldots,E_m$ such that for some cyclic order of $[n]$, such for every even $i$, $|E_{i-1}\cap E_i|=\ell$ and $|E_i\cap E_{i+1}=k-\ell$. Every $\ell$-offset Hamilton cycle consists of two perfect matchings of size $n/k$ and so $m=2n/k$. Dudek and Helenius proved: (i) if $k\geq 3$ and $1\leq \ell\leq k/2$ and $p\leq (1-\e)(e^k\ell!(k-1)!n^{-k})^{1/2}$ then w.h.p. $H_{n,p:k}$ does not contain an $\ell$-offset hamilton cycle; (ii) if $k\geq 3$ and $1\leq \ell\leq k/2$ and $p\geq (1+\e)(e^k\ell!(k-1)!n^{-k})^{1/2}$ then w.h.p. $H_{n,p:k}$ contains an $\ell$-offset hamilton cycle; (iii) if $k\geq 4$ and $\ell=2$ and $n^{k/2}p\to\infty$ then w.h.p. $H_{n,p:k}$ contains an 2-offset hamilton cycle.
\section{A related topic: long paths and cycles}
Erd\H{o}s conjectured that if $c>1$ then w.h.p. $G_{n,c/n}$ contains a path of length $f(c)n$ where $f(c)>0$. This was proved by Ajtai, Koml\'os and Szemer\'edi \cite{AKS1} and in a slightly weaker form by de la Vega \cite{Vega} who proved that if $c>4\log 2$ then $f(c)=1-O(c^{-1})$. See also Suen \cite{Suen}.
\begin{Problem}
Determine the precise form of $f(c)$ for $c$ close to one.
\end{Problem}
In particular, consider the case where $p=(1+\e)/n$ and where $\e^3n\to\infty$. {\L}uczak \cite{LucH} proved that w.h.p. the length of a longest path in $G_{n,p}$ lies in $(1+o(1))\e^2n[4/3,2]$. Anastos \cite{AShort} increased the lower bound to $\gtrsim 1.581\e^2n$.

Bollob\'as \cite{Bopath1} realised that for large $c$ one could find a large cycle w.h.p. by concentrating on a large subgraph with large minimum degree. In this way he showed that $f(c)\geq 1-e^{-24}c/2$. This was then improved by Bollob\'as, Fenner and Frieze \cite{BFF2} to $f(c)\geq 1-c^6e^{-c}$ and then by Frieze \cite{Fpath} to $f(c)\geq 1-(1+\e_c)(1+c)e^{-c}$ where $\e_c\to0$ as $c\to\infty$. This last result is optimal up to the value of $\e_c$, as there are w.h.p. $\approx (1+c)e^{-c}n$ vertices of degree 0 or 1. The paper \cite{Fpath} actually shows that for large $c$ there is w.h.p. a subgraph with property $\cA_k$ that contains most of the vertices of degree $k$ or more. Anastos and Frieze \cite{AFlongcycle} proved that the length of the longest cycle is w.h.p. $\approx f(c)n$ for some rather complicated definition of $f$. In particular $f(c)=1-(c+1)e^{-c}-c^2e^{-2c}+O(c^3e^{-3c})$. They also showed that  if $L_{c,n}$ is the length of the longest cycle in $G_{n,c/n}$, then $\E(L_{c,n}/n)$ tends to a limit. Using McDiarmid's coupling \cite{McD0} we see that a sparse random digraph has a path of length $\approx f(c)n$ w.h.p.
\begin{Problem}\label{Plong}
Explicitly determine $f(c)$ in the form $1-\sum_{k=1}^{\infty}p_k(c)e^{-c}$ where $p_k(c)$ is a polynomial in $c$.
\end{Problem}
Krivelevich, Lubetzky and Sudakov \cite{KLSB} proved that w.h.p. the random digraph $D_{n,c/n},c>1$ contains a cycle of length $n(1-(2+\e_c)e^{-c}$ vertices, where $\e_c\to0$ as $c\to\infty$. This is optimal up  to the value of $\e_c$. Anastos and Frieze \cite{AFlongdircycle} extended the result of \cite{AFlongcycle} to digraphs. They show that $D_{n,c/n},c>1$ contains a cycle of length $n\vec{f}$ where $\vec{f}(c)=1-2e^{-c}-(c^2+2c-1)e^{-2c}-O(c^3e^{-3c})$.
\begin{Problem}
Determine the precise form of $\vec{f}(c)$ for $c$ close to one.
\end{Problem}
\begin{Problem}
Explicitly determine $\vec{f}(c)$ in the form $1-\sum_{k=1}^{\infty}\vec{p}_k(c)e^{-c}$ where $\vec{p}_k(c)$ is a polynomial in $c$.
\end{Problem}
\begin{Problem}
Show that for large $c$, w.h.p. $D_{n,c/n}$ contains a subgraph containing most of the vertices with in- and out-degree $k$ and $k$ edge-disjoint Hamilton cycles.
\end{Problem}
Krivelevich, Kronenburg and Mond \cite{KKM} discuss the following Tur\'an question: given a random (or psedo-random) graph $G$ with $m$ edges and $t\in[n]$, what is the small value of $\a\in [0,1]$ such that every subgraph of $G$ with at least $\a m$ edges contains a cycle of length $t$?

Krivelevich, Lee and Sudakov \cite{KLS4} proved that if the graph $G$ has minimum degree $k$ and $kp\gg1$ then $G_p$ contains a cycle of length $(1-o(1))k$ with probability $1-o(1)$, (here $o(1)\to0$ as $k\to\infty$.) Furthermore, if $kp\geq (1+o(1))\log k$ then $G_p$ contains a cycle of length $k$ with probability $1-o(1)$. Riordan \cite{Riopath} gave a shorter proof of the first result of \cite{KLS4}.

If $H_{k-out}$ is as defined in Section \ref{kout}, then Frieze and Johansson \cite{FrJo} proved that if $H$ has minimum degree $m$ and $k$ is sufficiently large, then $H_{k-out}$ contains a cycle of length $(1-\e)m$ with probability $(1-o(1))$ (here $o(1)\to0$ as $m\to\infty$.)

Frieze and Jackson \cite{FrJa} considered the existence of large chordless cycles ({\em holes}) in the random graph $G_{n,c/n}$ and the random regular graph $G_{n,r}$. They proved that if $c$ is sufficiently large, then w.h.p. $G_{n,c/n}$ contains a hole of size $\Omega(n/c)$. They also proved that for every $r\geq 3$, $G_{n,r}$ contains a hole of size $\th_rn$ for some constant $\th_r>0$. Dragani\'c, Glock and Krivelevich \cite{DGK} showed that w.h.p. $G_{n,c/n}$ contains a hole of size $\approx\frac{2\log c}{c}n$, which is asymptotically optimal.
\begin{Problem}
Determine the size of the largest hole in $G_{n,r}$.
\end{Problem}
Dragani\'c, Glock and Krivelevich \cite{DGK1} give a short proof that w.h.p. $G_{n,p},p=(1+\e)/n$ contains an induced path of length $\Theta(\e^2n)$. Enriquez, Faraud, M\'enard and Noiry \cite{EFMN} gave an improved lower bound on the length of the longest induced cycle by giving an asymptotically sharp estimate of the length of the cycle constructed by the algorithm described in \cite{FrJa}.

Noever and Steger \cite{NoSt} showed that if $p=n^{-1/2+\e}$ then w.h.p. every subgraph of $G_{n,p}$ with minimum degree $(2/3+\e)np$ contains the square of a Hamilton cycle. \v{S}kori\'c, Steger and Truji\'c \cite{SST} improved this to show that if $p\geq C\bfrac{\log n}{n}^{1/k}, C=C(\a,\e)$, then w.h.p. every subgraph of $G_{n,p}$ with minimum degree at least $\brac{\frac{k}{k+1}+\a}np$ contains the $k$th power of a cycle on at least $(1-\e)n$ vertices.

Frieze \cite{Fcycles} showed that if $k|n,\,k=O(1)$ then w.h.p. a random graph withminimum degree at least two contains $k$ vertex disjoint cycles of size $n/k$ that cover $[n]$. The works of Johansson, Kahn, Vu \cite{JKV}, Kahn \cite{Ksham} and Hecke \cite{Heck} and Riordan \cite{Rio} give the threshold for the existence of a partition of $[n]$ into $n/3$ vertex disjoint triangles, assuming $3|n$.
\begin{Problem}
What is the threshold for $G_{n,p}$ to simultaneously contain for all partitions $k_1,k_2,\ldots,k_m$ of the integer $n$, vertex disjoint cycles $C_1,C_1,\ldots,C_m$ such that $|C_i|=k_i,i=1,2,\ldots,m$.
\end{Problem}
Alon, Krivelevich and Lubetzky \cite{AKL} study the set $\cL(G)$ of cycle lengths in sparse random graphs. They study random regular graphs as well as $G_{n,p}$. In the case of random regular graphs they establish the limiting probability that $\cL(G)\supset [\ell,n]$ for every $\ell\geq 3$. The results for $G_{n,p}$ are naturally slightly weaker, as the maximum of $\cL(G_{n,p})$ is more complicated. They show that w.h.p. $\cL(G_{n,p})\supseteq[\om(1),(1-\e)L_{c,n}]$ and the result of Anastos \cite{AA} shows that w.h.p. $\cL(G_{n,p})\supseteq[(1-\e)L_{c,n},L_{c,n}]$ completing the picture.
\begin{Problem}\label{last}
Discuss this problem in the context of other models of random graphs and hypergraphs.
\end{Problem}
Alon and Krivelevich \cite{AAKK} discuss the path cover number $\m(G)$ of $G=G_{n,p}$. By this we mean the minimum number of vertex disjoint paths needed to cover all vertices. They prove that for every $\e>0$ there exists $c_\e$ such if $c\geq c_\e$ and $p=c/n$ then w.h.p. $\m(G_{n,p})\in(1\pm\e)ce^{-c}n/2$. Alon and Anastos \cite{Alastos} sharpened the estimate of $\m(G_{n,p})$ quite considerably. They show 
that there is a function $f(c)$ such that $\m(G_{n,p=c/n})\approx f(c)n$ w.h.p. They show that $f(c)=\frac12ce^{-c}+e^{-c}+\cdots$. They also show that adding $\m$ edges will w.h.p. make the graph pancyclic.

Anastos, Erde, Kang and Pfenninger \cite{AEKP2025} proved that if $c\geq 20$ then there is a Central Limit Theorem for the length of the longest path in $G_{n,c/n}$.
\begin{Problem}
Remove the constraint $c\geq 20$.
\end{Problem}

Kozhevnikov, Raigorodskii and Zhukovskii \cite{KOZ} discussed the existence of long cycles in random subgraphs of Johnson graphs. Erde, Kang and Krivelevich \cite{EKK} proved that w.h.p. the random subgraph $Q_{n,p},p=c/n,\,c>1$ of the $n$-cube $Q_n$ contains a path of length $\Omega(2^n/(n\log n)^3)$. Anastos, Diskin, Erde, Kang, Krivelevich and Lichev \cite{ADEKKL} improved this to $2^n(1-o_c(1))$. In \cite{ADEKKLa}, the same set of authors proved that w.h.p. there are cycles of all lengths between 4 and $(1-\e)2^d$.
\begin{Problem}
Prove that $Q_{n,p},p=c/n$ contains a path of length $2^n(1-e^{-\Omega(c)})$ w.h.p
\end{Problem}
\section{Summary}
We have given a hopefully up to date description of what is known about Hamilton cycles in random graphs and hypergraphs. We have omitted extensions to pseudo-random graphs and other related topics. We have given 81 problems, some of which are a bit contrived. Here is a list of 9 which seem most interesting to me at the moment: \ref{P3}, \ref{P4}, \ref{P9}, \ref{P10}, \ref{P21}, \ref{P31}, \ref{P35}, \ref{P49}, \ref{Plong}. 


\begin{thebibliography}{999}
\bibitem{AigH} E. Aigner-Horev and D. Hefetz, \href{https://arxiv.org/pdf/2004.08637.pdf}{Rainbow hamilton cycles in randomly coloured randomly perturbed dense graphs}.
%
\bibitem{AKS1}  M. Ajtai, J. Koml\'{o}s and E. Szemer\'{e}di, The longest path in a random graph, {\em Combinatorica} 1 (1981) 1-12.
%
\bibitem{AKS}  M. Ajtai, J. Koml\'{o}s and E. Szemer\'{e}di, The first occurrence of Hamilton cycles in random graphs, {\em Annals of Discrete Mathematics} 27 (1985) 173-178.
%
\bibitem{ABKP} P. Allen,  J. B\"ottcher,  Y. Kohayakawa and Y. Person, Tight Hamilton cycles in random hypergraphs, {\em Random Structures and Algorithms} 45 (2015) 446-465.
%
\bibitem{AKPP} P. Allen,  C. Koch,  O. Parczyk and Y. Person, \href{https://arxiv.org/pdf/1710.08988.pdf}{Finding tight Hamilton cycles in random hypergraphs faster}.
%
\bibitem{Yalon} Y. Alon, \href{https://arxiv.org/pdf/2206.15235.pdf}{The global resilience of Hamiltonicity in $G(n,p)$}.
%
\bibitem{Alastos} Y. Alon and M. Anastos, \href{https://arxiv.org/pdf/2304.03710.pdf}{The completion numbers of Hamiltonicity and pancyclicity in random graphs}.
%
\bibitem{AK} Y. Alon and M. Krivelevich, \href{https://arxiv.org/pdf/1810.04987.pdf}{Random graph's Hamiltonicity is strongly tied to its minimum degree}.
%
\bibitem{AK1} Y. Alon and M. Krivelevich, \href{https://arxiv.org/pdf/1903.03007.pdf}{
Finding a Hamilton cycle fast on average using rotations-extensions}.
%
\bibitem{AKnew} Y. Alon and M. Krivelevich, \href{https://arxiv.org/pdf/1912.01251.pdf}{Hitting time of edge disjoint Hamilton cycles in random subgraph processes on dense base graphs}.
%
\bibitem{AAKK} Y. Alon and M. Krivelevich, \href{https://arxiv.org/pdf/2210.11770.pdf}{Hamilton completion and the path cover number of sparse random graphs}.
%
\bibitem{AKpan} Y. Alon and M. Krivelevich, \href{https://arxiv.org/pdf/2308.01564.pdf}{Sparse pancyclic subgraphs of random graphs}. 
%
\bibitem{AKpan1} Y. Alon and M. Krivelevich, \href{https://arxiv.org/pdf/2304.03710.pdf}{The completion numbers of Hamiltonicity and pancyclicity in random graphs}. 
%
\bibitem{AKL} Y. Alon, M. Krivelevich and E. Lubetzky, \href{https://arxiv.org/pdf/2008.13591.pdf}{Cycle lengths in sparse random graphs}.
%
\bibitem{AGIR} D. Altman, C. Greenhill, M. Isaev and R. Ramadurai, \href{https://arxiv.org/pdf/1611.09423.pdf}{A threshold result for loose Hamiltonicity in random regular uniform hypergraphs}.
%
\bibitem{AL}  A. Amit and N. Linial, Random Graph Coverings I: General Theory and Graph Connectivity, {\em Combinatorica} 22 (2002) 1-18.
%
\bibitem{AKLMS} J. Alvarado, Y. Kohayakawa, R. Lang, G. Mota and H. Stagni \href{https://arxiv.org/pdf/2309.14197.pdf}{RESILIENCE FOR LOOSE HAMILTON CYCLES}. 
%
\bibitem{AA} M. Anastos, \href{https://arxiv.org/pdf/2105.13828.pdf}{A note on long cycles in sparse random graphs}.
%
\bibitem{Ana22a} M. Anastos, \href{https://arxiv.org/pdf/2111.14771.pdf}{A fast algorithm on average for solving the Hamilton cycle problem}.
%
\bibitem{Ana22b} M. Anastos, \href{https://arxiv.org/pdf/2111.14759.pdf}{Fast algorithms for solving the Hamilton cycle problem with high probability}.
%
\bibitem{Ana22c} M. Anastos, \href{https://arxiv.org/pdf/2107.03527.pdf}{Packing Hamilton cycles in cores of random graphs}
%
\bibitem{AAA} M. Anastos, \href{https://arxiv.org/pdf/2209.09860.pdf}{Constructing Hamilton cycles and perfect matchings efficiently}.
%
\bibitem{AShort} M. Anastos, \href{https://arxiv.org/pdf/2208.06851.pdf}{An improved lower bound on the length of the longest cycle in random graphs}.
%
\bibitem{AnBr} M. Anastos and J. Briggs, \href{https://arxiv.org/pdf/1612.01965.pdf}{Packing Directed and Hamilton Cycles Online}.
%
\bibitem{AnCh} M. Anastos and D. Chakraborti, \href{https://arxiv.org/pdf/2309.12607.pdf}{ROBUST HAMILTONICITY IN FAMILIES OF DIRAC GRAPHS}. 
%
\bibitem{ADEKKL} M. Anastos, S. Diskin, J. Erde, M. Kang, M. Krivelevich and L. Lichev, \href{https://arxiv.org/pdf/2505.04436}{Nearly spanning cycle in the percolated hypercube}
%
\bibitem{ADEKKLa} M. Anastos, S. Diskin, J. Erde, M. Kang, M. Krivelevich and L. Lichev, \href{https://arxiv.org/pdf/2506.16858}{Cycle lengths in the percolated hypercube}.
%
\bibitem{AF} M. Anastos and A.M. Frieze, \href{https://arxiv.org/abs/1709.09198}{Pattern Colored Hamilton Cycles in Random Graphs}.
%
\bibitem{AEKP2025} M. Anastos, J. Erde, M. Kang and V. Pfenninger, \href{https://arxiv.org/pdf/2503.14336}{The law of the circumference of sparse binomial random graphs}.
%
\bibitem{AF1} M. Anastos and A.M. Frieze, \href{https://arxiv.org/pdf/1802.00433.pdf}{How many randomly colored edges make a randomly colored dense graph rainbow hamiltonian or rainbow connected?}.
%
\bibitem{AF19} M. Anastos and A.M. Frieze, \href{https://arxiv.org/pdf/1906.01433.pdf}{Hamilton cycles in random graphs with minimum degree at least 3: an improved analysis}.
%
\bibitem{AFlongcycle}  M. Anastos and A.M. Frieze, \href{https://arxiv.org/pdf/1907.03657.pdf}{A scaling limit for the length of the longest cycle in a sparse random graph}.
%
\bibitem{AFlongdircycle}  M. Anastos and A.M. Frieze, A scaling limit for the length of the longest cycle in a sparse random digraph.
%
\bibitem{AFG} M. Anastos, A.M. Frieze and J. Gao, \href{https://arxiv.org/pdf/1910.12594.pdf}{Hamiltonicity of random graphs in the stochastic block model}. 
%
\bibitem{AV} D. Angluin and L. Valiant, Fast probabilistic algorithms for hamiltonian circuits and matchings, {\em Journal of Computer and System Sciences} 18 (1979) 155-193.
%
\bibitem{ADS} S. Antoniuk, A. Dudek and A. Ruciński, \href{https://arxiv.org/pdf/2204.10738.pdf}{Powers of Hamiltonian cycles in randomly augmented Dirac graphs -- the complete collection}
%
\bibitem{ADS1} S. Antoniuk, A. Dudek and A. Ruciński, \href{https://arxiv.org/pdf/2512.23886}{Powers of Hamiltonian cycles in randomly augmented P\'osa-Seymour graphs}.
%
\bibitem{ADRRS} S. Antoniuk, A. Dudek, C. Reiher, A. Ruciński and M. Schacht, \href{https://arxiv.org/pdf/2002.05816.pdf}{High powers of Hamiltonian cycles in randomly augmented graphs}.
%
\bibitem{ABKPT} I. Araujo, J. Balogh, R. Krueger, S. Piga and A. Treglown, \href{https://arxiv.org/pdf/2212.10112.pdf}{On oriented cycles in randomly perturbed graphs}.
%
\bibitem{ABKPT1} I. Araujo, J. Balogh, R. Krueger, S. Piga and A. Treglown \href{https://arxiv.org/pdf/2212.10112.pdf}{ON ORIENTED CYCLES IN RANDOMLY PERTURBED DIGRAPHs}.
%
\bibitem{BBPP} D. Bal, P. Bennett, X. P\'erez-Gim\'enez and Pralat, Rainbow perfect matchings and Hamilton cycles in the random geometric graph, {\em Random Structures and Algorithms} 51 (2017) 587-606.
%
\bibitem{BalD} D. Bal and P. Devlin, \href{https://arxiv.org/pdf/1809.03596.pdf}{Hamiltonian Berge cycles in random hypergraphs}.
%
\bibitem{BBDS} D. Bal , R. Berkowitz, P. Devlin and M. Schacht,  \href{https://arxiv.org/pdf/1809.03596.pdf}{Hamiltonian Berge Cycles in Random Hypergraphs}.

%
\bibitem{BalF} D. Bal and A.M. Frieze, Rainbow Matchings and Hamilton Cycles in Random Graphs
{\em Random Structures and Algorithms} 48 (2016) 503-523.
%
\bibitem{BBSW} P. Balister, B. Bollob\'s, A. Sarkar and M. Walters. A critical constant for the $k$-nearest-neighbour model, {\em Advances in Applied Probability} 41 (2009) 1–12.
%
\bibitem{BBKMW} J. Balogh, B. Bollob\'as, M. Krivelevich, T. M\"uller and M. Walters, Hamilton cycles in random geometric graphs, {\em Annals of Applied Probability} 21 (2011) 1053-1072.
%
\bibitem{BHKM} W. Bedenknecht, J. Han, Y. Kohayakawa and  G. Mota, \href{https://arxiv.org/pdf/1802.08900.pdf}{Powers of tight Hamilton cycles in randomly perturbed hypergraphs}.
%
\bibitem{BellFrieze} T. Bell and A.M. Frieze, Rainbow powers of a Hamilton cycle in $G(n,p)$, {\em Journal of Graph Theory} 105 (2024) 491-500.
%
\bibitem{BFHK12} S. Ben-Shimon, A. Ferber, D. Hefetz and M. Krivelevich, Hitting time results for Maker-Breaker games, {\em Random Structures and Algorithms} 41 (2012) 23-46.
%
\bibitem{BKS}  S. Ben-Shimon, M. Krivelevich and B. Sudakov, On the resilience of Hamiltonicity and optimal packing of Hamilton cycles in random graphs, {\em SIAM Journal of Discrete Mathematics} 25 (2011) 1176-1193.
%
\bibitem{BKS2} S. Ben-Shimon, M. Krivelevich and B. Sudakov, Local resilience and Hamiltonicity Maker-Breaker games in random regular graphs, {\em Combinatorics, Probability and Computing} 20 (2011) 173-211.
%
\bibitem{BoF} T. Bohman and A.M. Frieze, Hamilton cycles in 3-out, {\em Random Structures and Algorithms} 35 (2009) 393-417.
%
\bibitem{BFM} T. Bohman, A.M. Frieze and R. Martin, How many random edges make a dense graph Hamiltonian? {\em Random Structures and Algorithms} 22 (2003) 33-42.
%
\bibitem{BR} M. Bloznelius and I. Radava\v{c}ius, A note on hamiltonicity of uniform random intersection graphs, {\em Lithuanian Mathematical Journal} 51 (2011).
%
\bibitem{Boll1} B. Bollob\'as, Random Graphs, First Edition, Academic Press, London 1985, Second Edition, Cambridge University Press, 2001.
%
\bibitem{Bopath1} B. Bollob\'as, Long paths in sparse random graphs, {\em Combinatorica} 2 (1982) 223-228.
%
\bibitem{Boll2}  B. Bollob\'as, The evolution of sparse graphs, Graph theory and combinatorics (Cambridge,  1983), Academic Press, London, (1984) 35–57.
%
\bibitem{BolRegHam} B. Bollob\'as, Almost all regular graphs are Hamiltonian, {\em European Journal of Combinatorics} 4 (1983) 97-106.
%
\bibitem{Bolmatch} B. Bollob\'as, Complete Matchings in Random Subgraphs of the Cube, {\em Random Structures and Algorithms} 1 (1990) 95-104.
%
\bibitem{BF} B. Bollob\'as and A. M. Frieze, On matchings and hamiltonian cycles in random graphs, {\em Annals of Discrete Mathematics} 28 (1985) 23-46.
%
\bibitem{BCFF} B. Bollob\'as, C. Cooper, T. Fenner and A.M. Frieze, On Hamilton cycles in sparse random graphs with minimum degree at least $k$, {\em  Journal of Graph Theory} 34 (2000) 42-59.
%
\bibitem{BFF} B. Bollob\'as, T. Fenner and A. M. Frieze, An algorithm for finding hamilton paths and cycles in random graphs, {\em Combinatorica} 7 (1987) 327-341.
%
\bibitem{BFF1} B. Bollob\'as, T. Fenner and A. M. Frieze, Hamilton cycles in random graphs with minimal degree at least $k$, {\em In A tribute to Paul Erd\H{o}s, Edited by A. Baker, B. Bollobás, A. Hajnal}, Cambridge University Press (1990) 59 - 96.
%
\bibitem{BFF2} B. Bollob\'as, T. Fenner and A. M. Frieze, Long cycles in sparse random graphs, {\em Graph theory and combinatorics, 59-64. Proceedings of Cambridge Combinatorial Conference in honour of Paul Erd\H{o}s} (1984) 59-64.
%
\bibitem{BoKo} B. Bollob\'as and Y. Kohayakawa, The hitting time of
Hamilton cycles in random bipartite graphs, {\em in Graph theory, combinatorics, algorithms, and applications (San Francisco, CA, 1989), SIAM} (1991) 26–41, 
%
\bibitem{BMPP} J. B\"ottcher, R. Montgomery, O. Parczyk and Yury Person, \href{https://arxiv.org/pdf/1802.04603.pdf}{Embedding Spanning Bounded Degree Graphs in Randomly Perturbed Graphs}.
%
\bibitem{BFKLS} J. Briggs, A.M. Frieze, M. Krivelevich, P. Loh abd B. Sudakov,  Packing Hamilton Cycles Online, {\em Combinatorics, Probability and Computing} 27 (2018) 475-496.  
%
\bibitem{BCCF} K. Burgin, P. Chebolu, C. Cooper and A.M. Frieze, Hamilton Cycles in Random Lifts of Graphs, {\em European Journal of Combinatorics} 27 (2006) 1282-1293.
%
\bibitem{CFH} D. Chakraborti, A.M. Frieze and M. Hasabanis, \href{https://arxiv.org/pdf/2103.03916.pdf}{Colorful Hamilton cycles in random graphs}.
%
\bibitem{CHS} Y. Chang, J. Han and L. Sun, \href{https://arxiv.org/pdf/2310.18980.pdf}{THE THRESHOLD FOR POWERS OF TIGHT HAMILTON CYCLES IN RANDOM HYPERGRAPHS}.
%
\bibitem{CHT}  Y. Chang, J. Han and L. Thoma, \href{https://arxiv.org/pdf/2007.11775.pdf}{On powers of tight Hamilton cycles in randomly perturbed hypergraphs}.
%
\bibitem{Chef} P. Chebolu and A.M. Frieze, Hamilton cycles in random lifts of complete directed graphs, {\em SIAM Journal on Discrete Mathematics} 22 (2008) 520-540.
%
\bibitem{CMM2025} M. Christoph, A. Martinsson and A. Milojevic, \href{https://arxiv.org/pdf/2505.05385}{Universality for transversal Hamilton cycles in random graphs}.
%
\bibitem{CNP} M. Christoph, R. Nenadov and K. Petrova, \href{https://arxiv.org/pdf/2402.01447}{The Hamilton space of pseudorandom graphs}.
%
\bibitem{CEP} D. Clemens, J. Ehrenm\"uller and Person, \href{https://arxiv.org/pdf/1903.09057.pdf}{A Dirac-type theorem for Berge cycles in random hypergraphs}.
%
\bibitem{CDGKO} P. Condon, A. Espuny D\'iaz, A. Gir$\tilde{\text{a}}$o, D. K\"uhn and  D. Osthus, \href{https://arxiv.org/pdf/1903.05052.pdf}{Dirac's theorem for random regular graphs}.
%
\bibitem{CDKKO1} P. Condon, A. Espuny D\'iaz, J. Kim, D. K\"uhn and  D. Osthus, \href{https://arxiv.org/pdf/1810.12433.pdf}{Resilient degree sequences with respect to Hamilton cycles and matchings in random graphs}.
%
\bibitem{CDGKO2} P. Condon, A. Espuny D\'iaz, Gir$\tilde{\text{a}}$o, D. K\"uhn and  D. Osthus, \href{https://arxiv.org/pdf/2007.02891.pdf}{Hamiltonicity of random subgraphs of the hypercube}.
%
\bibitem{Cook} N. Cook, L. Goldstein and T. Johnson, \href{https://arxiv.org/pdf/1510.06013.pdf}{Size biased couplings and the spectral gap for random regular graphs}, {\em Annals of Probability} 46 (2018) 72-125.  
%
\bibitem{C1} C. Cooper, Pancyclic Hamilton cycles in random graphs, {\em Discrete Mathematics} 91 (1991) 141-148.
%
\bibitem{C1a} C. Cooper, 1-Pancyclic Hamilton cycles in random graphs, {\em Random Structures and Algorithms} 3 (1992) 277-288.
%
\bibitem{Cperm} C. Cooper, The union of two random permutations does not have a directed Hamilton cycle, {\em Random Structures and Algorithms} 17 (2001) 95-98. 
%
\bibitem{CF1} C. Cooper and A.M. Frieze, On the number of hamilton cycles in a random graph, {\em Journal of Graph Theory} 13 (1989) 719-735.
%
\bibitem{CF3inout} C. Cooper and A.M. Frieze, Hamilton cycles in a class of random directed graphs, {\em Journal of Combinatorial Theory} B 62 (1994) 151-163.
%
\bibitem{CF-NN} C. Cooper and A.M. Frieze, On the connectivity of random $k$-th nearest neighbour graphs, {\em Combinatorics, Probability and Computing} 4 (1996) 343-362.
%
\bibitem{CF2inout} C. Cooper and A.M. Frieze, Hamilton cycles in random graphs and directed graphs, {\em Random Structures and Algorithms} 16 (2000) 369-401.
%
\bibitem{CF2} C. Cooper and A.M. Frieze, Multi-coloured Hamilton cycles in random edge-colored graphs, {\em Combinatorics, Probability and Computing} 11 (2002) 129-133.
 %
\bibitem{CF3} C. Cooper and A.M. Frieze, Multicoloured Hamilton cycles in random graphs: an anti-Ramsey threshold, {\em Electronic Journal of Combinatorics} 2 (1995).
%
\bibitem{CF4} C. Cooper and A.M. Frieze, Pancyclic random graphs, in {\em Proceedings of Random Graphs '87, Edited by M.Karonski, J.Jaworski and A.Rucinski} (1990) 29-39.
%
\bibitem{CoFr1} C. Cooper and A.M. Frieze, Hamilton cycles in random graphs and directed graphs, {\em Random Structures and Algorithms} 16 (2000) 368-401.
%
\bibitem{CFK} C. Cooper, A.M. Frieze and M. Krivelevich, Hamilton cycles in random graphs with a fixed degree sequence, {\em SIAM Journal on Discrete Mathematics} 24 (2010) 558-569.
%
\bibitem{CFM}  C. Cooper, A.M. Frieze and M. Molloy, Hamilton cycles in random regular digraphs,
{\em Combinatorics, Probability and Computing} 3 (1994) 39-50.
%
\bibitem{CFR} C. Cooper, A.M. Frieze and B. Reed, Random regular graphs of non-constant degree: connectivity and Hamilton cycles, {\em Combinatorics, Probability and Computing} 11 (2002) 249-262.
%
\bibitem{CFDham}  C. Cooper, and A.M. Frieze \href{https://arxiv.org/pdf/2312.06781.pdf}{Hamilton cycles in random digraphs with minimum degree at least one}.
%
\bibitem{DMP} J. D\'iaz, D. Mitsche and X. P\'erez-Gim\'enez, Sharp threshold for hamiltonicity of random geometric graphs, {\em SIAM Journal on Discrete Mathematics} 21 (2007) 57-65.
%
\bibitem{DGK} N. Dragani\'c, S. Glock and M. Krivelevich\href{https://arxiv.org/pdf/2106.00597.pdf}{The largest hole in sparse random graphs}.
%
\bibitem{DGK1} N. Dragani\'c, S. Glock and M. Krivelevich\href{https://arxiv.org/pdf/2106.08975.pdf}{Short proofs for long induced paths}.
%
\bibitem{DGCS} N. Dragani\'c, S. Glock, D. Corriea and B. Sudakov, \href{https://arxiv.org/pdf/2310.11580.pdf}{Optimal Hamilton covers and linear arboricity for random graphs}.
%
\bibitem{DK} N. Dragani\'c and P. Keevash, \href{https://arxiv.org/pdf/2502.00489}{P\'osa rotation through a random permutation}.
%
\bibitem{DEF}  A. Dudek, S. English and A.M. Frieze, On rainbow Hamilton cycles in random hypergraphs, {\em Electronic Journal of Combinatorics} 25 (2018).
%
\bibitem{DuFr1} A. Dudek and A.M. Frieze, Loose Hamilton Cycles in Random k-Uniform Hypergraphs,
{\em Electronic Journal of Combinatorics} 18, 2011,  
%
\bibitem{DuFr2} A. Dudek and A.M. Frieze, Tight Hamilton Cycles in Random Uniform Hypergraphs
{\em Random structures and Algorithms} 42 (2013) 374-385.
%
\bibitem{DFLS} A. Dudek,  A.M. Frieze, P. Loh and S. Speiss, Optimal divisibility conditions for loose Hamilton cycles in random hypergraphs, {\em Electronic Journal of Combinatorics} 19, 2012.
%
\bibitem{DFRS} A. Dudek, A.M. Frieze, A. Ruci\'nski and M. \v{S}ileikis, Embedding the Erd\H{o}s-R\'enyi Hypergraph into the Random Regular Hypergraph and Hamiltonicity, {Journal of Combinatorial Theorey }B (2017) 719-740.
%
\bibitem{DuH} A. Dudek and L. Helenius, On offset Hamilton cycles in random hypergraphs, {\em Discrete Applied Mathematics} 238 (2018) 77-85.
%
\bibitem{DRRS}  A. Dudek, C. Reiher, A. Ruci\'nski and M. Schacht, \href{https://arxiv.org/pdf/1805.10676.pdf}{Powers of Hamiltonian cycles in randomly augmented graphs}.
%
\bibitem{ES} C. Efthymiou and P. Spirakis, On the existence of Hamiltonian cycles in random intersection graphs, {\em in Proceedings of ICALP 32} (2005) 690-701. 
%
\bibitem{EFMN} N. Enriquez, G. Faraud, L. M\'enard and N. Noiry, \href{https://arxiv.org/pdf/2106.11130.pdf}{Long induced paths in a configuration model}.
%
\bibitem{EKK} J. Erde, M. Kang and M. Krivelevich, \href{https://arxiv.org/pdf/2106.04249.pdf}{Expansion, long cycles, and complete minors in supercritical random subgraphs of the hypercube}.
%
\bibitem{ER1} P. Erd\H{o}s and A. R\'enyi, On random graphs I, {\em Publ. Math. Debrecen} 6 (1959) 290-297. (1960) 17-61.
%
\bibitem{ER2} P. Erd\H{o}s and A. R\'enyi, On the evolution of random graphs, {\em Publ. Math. Inst. Hungar. Acad. Sci.} 5 (1960) 17-61. 
%
\bibitem{ESS} P. Erd\H{o}s, M. Simonovits and V. S\'os, Anti-Ramsey Theorems, {\em Colloquia Mathematica Societatis J\'anos Bolyai} 10, Infinite and Finite Sets, Keszethley, 1973. 
%
\bibitem{EFK} L. Espig, A.M. Frieze and  M. Krivelevich, Elegantly colored paths and cycles in edge colored random graphs, {\em SIAM Journal on Discrete Mathematics} 32 (2018) 1585-1618.
%
\bibitem{ED} A. Espuny D\'iaz, \href{https://arxiv.org/pdf/2102.02321.pdf}{HAMILTONICITY OF GRAPHS PERTURBED BY A RANDOM GEOMETRIC GRAPH}
%
\bibitem{EH} A. Espuny D\'iaz and J. Hyde, \href{https://arxiv.org/pdf/2205.08971.pdf}{Powers of Hamilton cycles in dense graphs perturbed by a random geometric graph}
%
\bibitem{DG}  A. Espuny D\'iaz and A. Gir$\tilde{\text{a}}$o, \href{https://arxiv.org/pdf/2101.06689.pdf}{HAMILTONICITY OF GRAPHS PERTURBED BY A RANDOM REGULAR GRAPH}.
%
\bibitem{DJKO} A. Espuny D\'iaz, F. Joos, D. K\"uhn and D. Osthus, \href{https://arxiv.org/pdf/1803.09223.pdf}{ 
Edge correlations in random regular hypergraphs and applications to subgraph testing}.
%
\bibitem{DR}  A. Espuny D\'iaz and R. Razafindravola, \href{https://arxiv.org/pdf/2410.14447}{How many random edges make an almost-Dirac graph Hamiltonian?}.
%
\bibitem{FeFr1} T. Fenner and A.M. Frieze, On the existence of hamiltonian cycles in a class of random graphs, {\em Discrete Mathematics} 45 (1983) 301-305.
%
\bibitem{FFregHam} T. Fenner and A.M. Frieze, Hamiltonian cycles in random regular graphs, {\em Journal of Combinatorial Theory B} 37 (1984) 103-112.
%
\bibitem{FGKN15} A. Ferber, R. Glebov, M. Krivelevich and A. Naor, Biased games on random boards, {\em Random Structures and Algorithms}, 46 (2015) 651-676.  
%
\bibitem{FeKr} A. Ferber and  M. Krivelevich, Rainbow Hamilton cycles in random graphs and hypergraphs, in {\em Recent trends in combinatorics, IMA Volumes in Mathematics and its applications, A. Beveridge, J. R. Griggs, L. Hogben, G. Musiker and P. Tetali, Eds.}, Springer 2016, 167-189.
%
\bibitem{FKSV} A. Ferber, M. Krivelevich, B. Sudakov and P. Vieira, Finding paths in sparse random graphs requires many queries, {\em Random Structures and Algorithms} 49 (2016) 635-668.
%
\bibitem{FKL} A. Ferber, G. Kronenberg and E. Long, Packing, Counting and Covering Hamilton cycles in random directed graphs, {\em Israel Journal of Mathematics} 220 (2017) 57–87.
%
\bibitem{FKS} A. Ferber, M. Kwan and B. Sudakov, \href{https://arxiv.org/pdf/1708.07746.pdf}{Counting Hamilton cycles in sparse random directed graphs}.
%
\bibitem{FHM} A. Ferber, J. Han and D. Mao, \href{https://arxiv.org/pdf/2211.05477.pdf}{Dirac-type Problem of Rainbow matchings and Hamilton cycles in Random Graphs}.
%
\bibitem{FeLo} A. Ferber and E. Long, \href{https://arxiv.org/pdf/1603.03614.pdf}{Packing and counting arbitrary Hamilton cycles in random digraphs}.
%
\bibitem{FNNPS} A. Ferber, R. Nenadov, A. Noever, U. Peter and N. \v{S}kori\'c, Robust hamiltonicity of random directed graphs, {\em Journal of Combinatorial Theory} B 126 (2017) 1-23. 
%
\bibitem{FSS} A. Ferber, M. Sales and M. Shurman \href{https://arxiv.org/pdf/2410.12964}{Covering random digraphs with Hamilton cycles}.
%
\bibitem{FSST}  M. Fischer, N. \v{S}kori\'c, A. Steger and M. Truji\'c, \href{https://arxiv.org/pdf/1809.07534.pdf}{Triangle resilience of the square of a Hamilton cycle in random graphs}.
%
\bibitem{FMMS} N. Fountoulakis, D. Mitsche, T. M\"uller and M. Schepers, \href{https://arxiv.org/pdf/1901.09175.pdf}{Hamilton cycles and perfect matchings in the KPKVB model.}
%
\bibitem{FKNP} K. Frankston, J. Kahn, B. Narayanan and J. Park, \href{https://arxiv.org/pdf/1910.13433.pdf}{Thresholds versus fractional expectation-thresholds}
%
\bibitem{Fcycles} A.M. Frieze, Partitioning random graphs into large cycles, {\em Discrete Mathematics} 70 (1988) 149-158.
%
\bibitem{F0} A.M. Frieze, \href{http://www.math.cmu.edu/~af1p/Texfiles/REGNONCONSTANT.pdf}{Random regular graphs of non-constant degree.}
%
\bibitem{Fbip} A.M. Frieze, Limit distribution for the existence of hamiltonian cycles in random bipartite graphs, {\em European Journal of Combinatorics} 6 (1985) 327-334.
%
\bibitem{Fpath}  A.M. Frieze, On large matchings and cycles in sparse random graphs, {\em Discrete Mathematics} 59 (1986) 243-256.
%
\bibitem{Fdir} A.M. Frieze, An algorithm for finding hamilton cycles in random digraphs, {\em Journal of Algorithms} 9 (1988) 181-204.
%
\bibitem{F2} A.M. Frieze, Parallel algorithms for finding Hamilton cycles in random graphs, {\em Information Processing Letters} 25 (1987) 111-117.
%
\bibitem{Freg} A.M. Frieze, Finding hamilton cycles in sparse random graphs, {\em Journal of Combinatorial Theory B}  44 (1988) 230-250.
%
\bibitem{Fperm} A.M. Frieze, Hamilton cycles in the union of random permutations, {\em Random Structures and Algorithms} 18 (2001) 83-94.
%
\bibitem{F1} A.M. Frieze, On a Greedy 2-Matching Algorithm and Hamilton Cycles in Random Graphs with Minimum Degree at Least Three, {\em Random structures and Algorithms} 45 (2014) 443-497.
%
\bibitem{Fhyp1} A.M. Frieze, Loose Hamilton Cycles in Random 3-Uniform Hypergraphs, {\em Electronic Journal of Combinatorics} 17, 2010.
%
\bibitem{FH} A.M. Frieze and S. Haber, An almost linear time algorithm for finding Hamilton cycles in sparse random graphs with minimum degree at least three, {\em Random Structures and Algorithms} 47 (2015) 73-98.
%
\bibitem{FrJa} A.M. Frieze and B. Jackson, Large holes in sparse random graphs, {\em Combinatorica} 7 (1987) 265-274.
%
\bibitem{FJMRW} A.M. Frieze, M.R. Jerrum, M. Molloy, R. Robinson and N.C. Wormald, Generating and counting Hamilton cycles in random regular graphs, {\em Journal of Algorithms} 21 (1996) 176-198.
%
\bibitem{FrJo} A.M. Frieze and T. Johansson,  On random $k$-out sub-graphs of large graphs, {\em Random Structures and Algorithms} 50 (2017) 143-157.
%
\bibitem{FK} A.M. Frieze and M. Karo\'nski, Introduction to Random Graphs, Cambridge University Press,  2015.
%
\bibitem{FKT} A.M. Frieze, M. Karo\'nski and L. Thoma, On Perfect Matchings and Hamiltonian Cycles in Sums of Random Trees, {\em SIAM Journal on Discrete Mathematics} 12 (1999) 208-216.
%
\bibitem{FrKr}  A.M. Frieze and M. Krivelevich, On two Hamilton cycle problems in random graphs, {\em Israel Journal of Mathematics} 166 (2008) 221-234.
%
\bibitem{FKMP}  A.M. Frieze, M. Krivelevich, P. Michaeli and R. Peled, On the trace of random walks on random graphs, {\em Proceedings of the London Mathematical Society} 16 (2018) 847-877. 
%
\bibitem{builder} A.M. Frieze, M. Krivelevich and P. Michaeli, Fast construction on a restricted budget.
%
\bibitem{FLo} A.M. Frieze and P. Loh, Rainbow hamilton cycles in random graphs, {\em Random Structures and Algorithms} 44 (2014) 328-354.
%
\bibitem{FL} A.M. Frieze and T. {\L}uczak, Hamiltonian cycles in a class of random graphs: one step further, {\em n Proceedings of Random Graphs '87, Edited by M.Karonski, J.Jaworski and A.Rucinski, John Wiley and Sons} (1989) 53-59.
%
\bibitem{FM} A.M. Frieze and B. Mckay, Multicoloured trees in random graphs, {\em Random Structures and Algorithms} 5 (1994) 45-56.
%
\bibitem{FP} A.M. Frieze and X. P\'erez-Gim\'enez, \href{https://arxiv.org/pdf/2003.02998.pdf}{Rainbow Hamilton Cycles in Random Geometric Graphs}.
%
\bibitem{FP1} A.M. Frieze and X. P\'erez-Gim\'enez
%
\bibitem{FPgame} A.M. Frieze and W. Pegden,  \href{https://arxiv.org/pdf/2003.13521.pdf}{Maker Breaker on Digraphs}.
%
\bibitem{arbpattern}  A.M. Frieze and W. Pegden, \href{https://arxiv.org/pdf/2305.00880.pdf}{Sequentially constrained Hamilton Cycles in random graphs}.
%
\bibitem{FPPR} A.M. Frieze, P. Pralat, X. P\'erez-Gim\'enez and B. Reiniger, \href{https://arxiv.org/pdf/1610.07988.pdf}{Perfect matchings and Hamiltonian cycles in the preferential attachment model}.
%
\bibitem{FPP} A.M. Frieze, P. Pralat and X. P\'erez-Gim\'enez, \href{https://arxiv.org/pdf/2002.07313.pdf}{On the existence of Hamilton cycles with a periodic pattern in a random digraph}.
%
\bibitem{FGKMPS} A.M. Frieze, P. Gao, C. MacRury, P. Pralat and G. Sorkin,\href{https://arxiv.org/pdf/2311.05533}{Building Hamiltonian Cycles in the Semi-Random Graph Process in Less Than $2n$ Rounds}. \href{https://arxiv.org/pdf/2208.00255.pdf}{Hamilton cycles in a semi-random graph model}.
%
\bibitem{Ganeson} G. Ganeson, \href{https://arxiv.org/pdf/2112.05641.pdf}{Bridged Hamilton Cycles in Sub-Critical Random Geometric graphs}.
%
\bibitem{GIM} J. Gao, M. Isaev and B. McKay, \href{https://arxiv.org/pdf/1906.02886.pdf}{Sandwiching random regular graphs between binomial random graphs}
%
\bibitem{GKMP} P. Gao, B. Kaminski, C. MacRury and  P. Pralat, \href{https://arxiv.org/pdf/2006.02599.pdf}{Hamilton Cycles in the Semi-random Graph Process}.
%
\bibitem{GMP} P. Gao, C. MacRury and  P. Pralat, \href{https://arxiv.org/pdf/2205.02350.pdf}{A Fully Adaptive Strategy for Hamiltonian Cycles in the Semi-Random Graph Process}
%
\bibitem{GKM} L. Gishboliner, M. Krivelevich and P. Michaeli, \href{https://arxiv.org/pdf/2007.12111.pdf}{Colour-biased Hamilton cycles in random graphs}.
%
\bibitem{GKM1} L. Gishboliner, M. Krivelevich and P. Michaeli, \href{https://arxiv.org/pdf/2203.07148.pdf}{Oriented discrepancy of Hamilton cycles}
%
\bibitem{GK} R. Glebov and M. Krivelevich, On the number of Hamilton cycles in sparse random graphs, {\em SIAM Journal on Discrete Mathematics} 27 (2013) 27-42.
%
\bibitem{GKS} R. Glebov, M. Krivelevich and T. Szabo, On covering expander graphs by Hamilton cycles, {\em Random Structures and Algorithms} 44 (2014) 183-200.  
%
\bibitem{GNS} R. Glebov H. Naves and B. Sudakov, The Threshold Probability for Long Cycles, {\em Combinatorics, Probability and Computing} 26 (2017) 208-247.
%
\bibitem{GlSg} S. Glock and A. Sguegla, \href{https://arxiv.org/pdf/2510.01949}{On Kotzig's conjecture in random graphs}.
%
\bibitem{GS} Y. Gurevich and S. Shelah, Expected computation time for Hamiltonian path problem, {\em SIAM Journal on Computing} 16 (1987) 486-502. 
%
\bibitem{Heck} A. Heckel, \href{https://arxiv.org/pdf/1802.08472.pdf}{Random triangles in random graphs}.
%
\bibitem{HKLO} D. Hefetz, D. K\"{u}hn, J. Lapinskas and D. Osth\"{u}s, Optimal covers with Hamilton cycles in random graphs, {\em Combinatorica} 34 (2014) 573-596.
%
\bibitem{HSS} D. Hefetz, A. Steger and B. Sudakov, Random directed graphs are robustly Hamiltonian, {\em Random Structures and Algorithms} 49 (2016) 345-362.
%
\bibitem{HefKriv} D. Hefetz and M. Krivelevich, \href{https://arxiv.org/pdf/2506.19731}{The Hamilton cycle space of random graphs}. 
%
\bibitem{HefKriv1} D. Hefetz and M. Krivelevich, \href{https://arxiv.org/pdf/2507.04488}{The Hamilton cycle space of random regular graphs and randomly perturbed graphs}.
%
\bibitem{HKT} D. Hefetz, M. Krivelevich and W. E. Tan, Waiter-Client and Client-Waiter Hamiltonicity games on random graphs, {\em European Journal of Combinatorics} 63 (2017) 26-43. 
%
\bibitem{HK} M. Held and R. Karp, A Dynamic Programming Approach to Sequencing Problems, {\em SIAM Journal of Applied Mathematics} 10 (1962) 196-210.
%
\bibitem{HLMM} C. Henderson, S. Longbrake, D. Mao and P. Morawski, \href{https://arxiv.org/pdf/2506.21756}{Hamilton cycles in regular graphs perturbed by a random 2-factor}.
%
\bibitem{ImKim} S. Im and M. Kim, \href{https://arxiv.org/pdf/2512.06675}{Berge Hamilton cycles in a random sparsification of dense hypergraphs}.
%
\bibitem{JW} S. Janson and N. Wormald, Rainbow Hamilton cycles in random regular graphs, {\em 
Random Structures and Algorithms} 30 (2007) 35-49.
%
\bibitem{JS} M. Jerrum and A. Sinclair, Approximating the permanent, {\em SIAM Journal on Computing} 18 (1989) 1148-1178.
%
\bibitem{Ja1} S. Janson, The numbers of spanning trees, Hamilton cycles and perfect matchings in a random graph, {\em Combinatorics, Probability and Computing} 3 (1994) 97-126. 
%
\bibitem{JLR} S. Janson, T. {\L}uczak and A. Ruci\'nski, Random Graphs, John Wiley and Sons, New York, 2000.
%
\bibitem{JLS} F. Joos, R. Lang and N. Sanhueza-Matamala, \href{https://arxiv.org/pdf/2312.15262}{Robust Hamiltonicity}
%
\bibitem{KimVu} J. Kim and V. Vu, Sandwiching random graphs: universality between random graph models, {\em Advances in Mathematics} 188 (2004) 444-469.
%
\bibitem{KNP} J. Kahn, B. Narayanan and J. Park, \href{https://arxiv.org/pdf/2010.08592.pdf}{The threshold for the square of a Hamilton cycle}.
%
\bibitem{KW} J. Kim and N. Wormald, Random matchings which induce Hamilton cycles, and Hamiltonian decompositions of random regualr graphs, {\em Journal of Combinatorial Theory B} 81 (2001) 20-44.
%
\bibitem{JKV} A. Johansson, J. Kahn and V. Vu, Factors in random graphs, {\em Random Structures and Algorithms} 33 (2008) 1-28.
%
\bibitem{Jo} T. Johansson, On Hamilton cycles in Erd\H{o}s-R\'enyi subgraphs of large graphs, {\em Random Structures and Algorithms} 57 (2020) 132-149.
%
\bibitem{Jo1} T. Johansson, \href{https://arxiv.org/pdf/2001.05258.pdf}{A condition for Hamiltonicity in Sparse Random Graphs with a Fixed Degree Sequence}.
%
\bibitem{Ksham} J. Kahn, \href{https://arxiv.org/pdf/1909.06834.pdf}{Asymptotics for Shamir's Problem}.
%
\bibitem{J2020} T. Johansson, \href{https://arxiv.org/pdf/2012.11953.pdf}{Hamilton cycles in weighted Erd\H{o}s-R\'enyi graphs}
%
\bibitem{KaLe} K. Katsamaktsis and S. Letzter, \href{https://arxiv.org/pdf/2304.09155.pdf}{Rainbow Hamiltonicity in uniformly colored perturbed graphs}.
%
\bibitem{KaLS} K. Katsamaktsis, S. Letzter and A. Sgueglia,\href{https://arxiv.org/pdf/2310.18284.pdf}{Rainbow subgraphs of uniformly coloured randomly perturbed graphs}.
%
\bibitem{KaLS1} K. Katsamaktsis, S. Letzter and A. Sgueglia, \href{https://arxiv.org/pdf/2304.09155.pdf}{Rainbow Hamiltonicity in uniformly coloured perturbed digraphs}.
%
\bibitem{KMP}  Kelly, M\''uyesser and Pokrovskiy, \href{https://arxiv.org/pdf/2308.08535}{Optimal spread for spanning subgraphs of Dirac hypergraphs}.
%
\bibitem{KKO} F. Knox, D. K\"uhn and D. Osthus, Edge-disjoint Hamilton cycles in random graphs, {\em Random Structures and Algorithms} 46 (2015), 397-445.
%
\bibitem{KS1} J. Koml\'os and E. Szemer\'edi, Hamilton cycles in random graphs, {\em In: A. Hajnal, R. Rado, V. T. S\'os, eds., Infinite and Finite Sets, Colloq. Math. Soc. Janos Bolyai} 10 (North-Holiand, Amsterdam), 1973.
%
\bibitem{KS2} J. Koml\'os and E. Szemer\'edi, Limit distributions for the existence of Hamilton circuits in a random graph, {\em Discrete Mathematics} 43 (1983) 55-63.
%
\bibitem{Kor} A.D. Korsunov, Solution of a problem of Erd\H{o}s and R\'enyi on hamiltonian cycles in nonoriented graphs, Soviet Math. Doklaidy 17 (1976) 760-764.
%
\bibitem{KOZ} V. Kozhevnikov, A. Raigorodskii and M. Zhukovskii, \href{https://arxiv.org/pdf/2012.07314.pdf}{Large cycles in random generalized Johnson graphs}.
%
\bibitem{KKM} M. Krivelevich, G. Kronenberg and Adva Mond, \href{https://arxiv.org/pdf/1911.08539.pdf}{Tur\'an-type problems for long cycles in random and pseudo-random graphs}.
%
\bibitem{KKS} M. Krivelevich, M. Kwan and B. Sudakov, Cycles and matchings in randomly perturbed digraphs and hypergraphs, {\em Combinatorics, Probability and Computing} 25 (2016) 909-927.
%
\bibitem{KLS} M. Krivelevich, E. Lubetzky and B. Sudakov, Cores of random graphs are born Hamiltonian, {\em Proceedings of the London Mathematical Society} 109 (2014) 161-188.
%
\bibitem{KLS1} M. Krivelevich, C. Lee and B. Sudakov, Compatible Hamilton cycles in random graphs, {\em Random Structures and Algorithms} 49 (2016) 533-557.  
%
\bibitem{KLS2} M. Krivelevich, C. Lee and B. Sudakov, Resilient pancyclicity of random and pseudo-random graphs, {\em SIAM Journal on Discrete Mathematics} 24 (2010) 1-16. 
%
\bibitem{KLS3} M. Krivelevich, C. Lee and B. Sudakov, Robust Hamiltonicity of Dirac graphs, {\em Transactions of the American Mathematical Societ}y 366 (2014) 3095-3130.  
%
\bibitem{KLS4} M. Krivelevich, C. Lee and B. Sudakov, Long paths and cycles in random subgraphs of graphs with large minimum degree, {\em Random Structures and Algorithms} 46 (2015) 320-345.
%
\bibitem{KLSA} M. Krivelevich, E. Lubetzky and B. Sudakov, Hamiltonicity thresholds in Achlioptas processes., {\em Random Structures and Algorithms} 37 (2010) 1-24.  
%
\bibitem{KLSB}  M. Krivelevich, C. Lee and B. Sudakov,  Longest cycles in sparse random digraphs., {\em Random Structures and Algorithms} 43 (2013) 1-15.  
%
\bibitem{KSa} M. Krivelevich and W. Samotij, Optimal packings of Hamilton cycles in sparse random graphs, {\em SIAM Journal on Discrete Mathematics} 26 (2012) 964-982.
%
\bibitem{KSVW} M. Krivelevich, B. Sudakov, V. Vu and N. Wormald, Random regular graphs of high degree, {\em  Random Structures and Algorithms} 18 (2001) 346-363.
%
\bibitem{KO} D. K\"uhn and D. Osthus, On P\'osa's conjecture for random graphs, {\em SIAM Journal on Discrete Mathematics} 26 (2012) 1440-1457. 
%
\bibitem{KO1}  D. K\"uhn and D. Osthus, Hamilton decompositions of regular expanders: applications, {\em Journal of Combinatorial Theory  B} 104 (2014) 1-27. 
%
\bibitem{LS} C. Lee and B. Sudakov, Dirac's theorem for random graphs, {\em Random Structures and Algorithms} 41 (2012) 293-305.
%
\bibitem{LS1} C. Lee and W. Samotij, Pancyclic subgraphs of random graphs, {\em Journal of Graph Theory} 71 (2012) 142-158.
%
\bibitem{LLP} E. Levy, G. Louchard and J. Petit,  A Distributed Algorithm to Find Hamiltonian Cycles in $G_{n,p}$ Random Graphs, {\em Combinatorial and Algorithmic Aspects of Networking, Eds. A. L{\'o}pez-Ortiz and A. Hamel} (2005) 63-74.
%
\bibitem{LWW} T. {\L}uczak, {\L}. Witkowski and M. Witkowski, Hamilton cycles in random lifts of graphs, {\em European Journal of Combinatorics} 49 (2015) 105-116.
%
\bibitem{LucH} T. {\L}uczak, Cyles in a random graph near the critical point, {\em Random Structures and Algorithms} 2 (1991) 421-439.
%
\bibitem{MY} G. Ma and Z. Yan, \href{https://arxiv.org/pdf/2605.29553}{Sharp threshold for hamilton cycles in randomly perturbed sparse graphs}.
%
\bibitem{MPPS} T. Makai, M. Pasch, K. Petrova and L. Schiller, \href{https://arxiv.org/pdf/2502.14515}{Sharp thresholds for higher powers of Hamilton cycles in random graphs}.
%
\bibitem{MPW} T. M\"uller, X. P\'erez and N. Wormald, Disjoint Hamilton cycles in the random geometric graph, {\em Journal of Graph Theory} 68 (2011) 299-322.
%
\bibitem{McD0} C. McDiarmid, Clutter percolation and random graphs, {\em Mathematical Programming Studies} 13 (1980) 17-25.
%
\bibitem{McD1}  C. McDiarmid, Expected numbers at hitting times, {\em Journal of Graph Theory} 15 (1991) 637-648.
%
\bibitem{McM} A. McDowell and R. Mycroft, \href{https://arxiv.org/pdf/1802.04242.pdf}{Hamilton $\ell$-cycles in randomly-perturbed hypergraphs}.
%
\bibitem{MS} 	P. MacKenzie and Q. Stout, Optimal parallel construction of Hamiltonian cycles and spanning trees in random graphs, {\em  Proceedings of the fifth annual ACM symposium on Parallel algorithms and architectures} (1993) 224-229.
%
\bibitem{MPS} M. Molloy, P. Pralat and G. Sorkin, \href{https://arxiv.org/pdf/2401.00559}{Perfect matchings and loose Hamilton cycles in the semirandom hypergraph model}.
%
\bibitem{M1} R. Montgomery, \href{https://arxiv.org/abs/1710.00505}{Hamiltonicity in random graphs is born resilient}.
%
\bibitem{Mnotes}  R. Montgomery, \href{http://web.mat.bham.ac.uk/R.H.Montgomery/topicsinrandomgraphs.pdf}{Topics in random graphs}.
%
\bibitem{M2} R. Montgomery, \href{https://arxiv.org/pdf/1901.09605.pdf}{Hamiltonicity in random directed graphs is born resilient.}
%
\bibitem{M3} R. Montgomery, \href{https://arxiv.org/pdf/2103.06751.pdf}{Spanning cycles in random directed graphs}
%
\bibitem{NaSc} B. Narayanan and M. Schacht, \href{https://www.math.uni-hamburg.de/home/schacht/preprints/nonlin_thresh.pdf}{Sharp thresholds for nonlinear Hamiltonian cycles in hypergraphs}.
%
\bibitem{NS} R. Nenadov and N. \v{S}kori\'c, \href{https://arxiv.org/abs/1601.04034}{Powers of Hamilton cycles in random graphs and tight Hamilton cycles in random hypergraphs}.
%
\bibitem{NSU} R. Nenadov, A. Steger and P.Su, \href{https://arxiv.org/pdf/2012.02551.pdf}{An $O(n)$ time algorithm for finding Hamilton cycles with high probability}.
%
\bibitem{NST} R. Nenadov, A. Steger and M. Truji\'c \href{https://arxiv.org/abs/1710.00799}{Resilience of Perfect Matchings and Hamiltonicity in Random Graph Processes}.
%
\bibitem{NeTr}  R. Nenadov and M. Truji\'c, \href{https://arxiv.org/pdf/1811.09209.pdf}{Sprinkling a few random edges doubles the power}.
%
\bibitem{NoSt} A. Noever and A. Steger, \href{https://www.combinatorics.org/ojs/index.php/eljc/article/view/v24i4p8}{Local Resilience for Squares of Almost Spanning Cycles in Sparse Random Graphs}.
%
\bibitem{PaPe} O. Parczyk and Y. Person, \href{https://arxiv.org/pdf/1504.02243.pdf}{Spanning structures and universality in sparse hypergraphs}.
%
\bibitem{Pen} M. Penrose, Random Geometric Graphs, Oxford University Press, 2003.
%
\bibitem{PeTr} K. Petrova and M. Truji\'c, \href{https://arxiv.org/pdf/2205.11421.pdf}{Transference for loose Hamilton cycles in random 3-uniform hypergraphs}
%
\bibitem{Poole} D. Poole, \href{https://arxiv.org/pdf/1410.7446.pdf}{On Weak Hamiltonicity of a Random Hypergraph}.
%
\bibitem{Posa} L. P\'osa, Hamiltonian circuits in random graphs, {\em Discrete Mathematics} 14 (1976) 359-364.
%
\bibitem{R} O. Riordan, Spanning subgraphs of random graphs, {\em Combinatorics, Probability and Computing} 9 (2000) 125-148.
%
\bibitem{Riopath} O. Riordan, Long cycles in random subgraphs of graphs with large minimum degree, {\em  Random Structures Algorithms} 45 (2014) 762-765.
%
\bibitem{Rio} O. Riordan, \href{https://arxiv.org/pdf/1802.01948.pdf}{Random cliques in random graphs.}
%
\bibitem{RW} R. Robinson and N. Wormald, Existence of long cycles in random cubic graphs, in {\em Enumeration and Design, D. Jackson and S. Vaustone Eds.} (1984) 251-270.
%
\bibitem{RW1} R. Robinson and N. Wormald, Almost all cubic graphs are Hamiltonian, {\em Random Structures and Algorithms} 3 (1992) 117-125.
%
\bibitem{RWmat} R. Robinson and N. Wormald, Hamilton cycles containing randomly selected edges in random regular graphs, {\em Random Structures and Algorithms} 19 (2001) 128 - 147.
%
\bibitem{RW1a} R. Robinson and N. Wormald, Almost all regular graphs are Hamiltonian, {\em Random Structures and Algorithms} 5 (1994) 363-374.
%
\bibitem{RW2} R. Robinson and N. Wormald, Hamilton cycles containing randomly selected edges in random regular graphs, {\em Random Structures and Algorithms} 19 (2001) 128-147.
%
\bibitem{Ryb1} K. Rybarczyk, Sharp threshold functions for random intersection graphs via a coupling method, {\em The Electronic Journal of Combinatorics} 18 (2011).
%
\bibitem{Ryb2} K. Rybarczyk, \href{https://arxiv.org/pdf/1702.03667.pdf}{Finding Hamilton cycles in random intersection graphs}.
%
\bibitem{S} E. Shamir, How many random edges make a graph hamiltonian?, {\em Combinatorica} 3 (1983) 123–131.
%
\bibitem{SST} N. \v{S}kori\'c, A. Steger and M. Truji\'c, Local resilience of an almost spanning $k$-cycle in random graphs, {\em Random Structures and Algorithms} 53 (2018) 728-751.
%
\bibitem{ST} D. Spielman and S-H. Teng, Smoothed Analysis of Algorithms: Why The Simplex Algorithm Usually Takes Polynomial Time, {\em Journal of the ACM} 51 (2004) 385-463.
%
\bibitem{SV} B. Sudakov and V. Vu,  Local resilience of graphs, {\em Random Structures and Algorithms} 33 (2008) 409-433.
%
\bibitem{Suen} S. Suen, On large induced trees and long induced paths in sparse random graphs, {\em Journal of Combinatorial Theory B} 56 (1992) 250-262.
%
\bibitem{TWZ} A. Telcs, N. Wormald and S. Zhou, Hamiltonicity of random graphs produced by 2-processes, {\em Random Structures and Algorithms} 31 (2007) 450-481. 
%
\bibitem{T} A. Thomason, A simple linear expected time algorithm for finding a hamilton path, {\em Discrete Mathematics} 75 (1989) 373-379.
%
\bibitem{Tik} K. Tikhomirov and P. Youssef, The spectral gap of dense random regular graphs, {\em Annals of Probability} 47 (2019) 362-419.
%
\bibitem{Tu} V. Tureau, \href{https://arxiv.org/pdf/1805.06728.pdf}{A Distributed Algorithm for Finding Hamiltonian Cycles in Random Graphs in O(log n) Time}.
%
\bibitem{Vega} W. de la Vega, Long paths in random graphs, {\em Studia Scient. Math. Hungar.} 14 (1979) 335-340.
%
\bibitem{Wang} Y. Wang, \href{https://arxiv.org/pdf/2606.11992}{On the hitting time of hamiltonicity in bipartite Dirac graphs}.
%
\bibitem{Wormreg} N. Wormald, Models of random regular graphs, {\em in Surveys in Combinatorics, Edited by J. D. Lamb and D. A. Preece} (1999) 239-298.
%
\bibitem{Zuk} M. Zukoovskii, \href{https://arxiv.org/pdf/2502.14794}{Sharp thresholds for spanning regular subgraphs}.
\end{thebibliography}
\end{document}